\documentclass{article}
\usepackage{arxiv}

\usepackage[utf8]{inputenc} 
\usepackage[T1]{fontenc}    
\usepackage{hyperref}       
\usepackage{url}            
\usepackage{booktabs}       
\usepackage{amsfonts}       
\usepackage{nicefrac}       
\usepackage{microtype}      
\usepackage{lipsum}
\usepackage{graphicx}

\usepackage{algorithm}
\usepackage{algorithmicx}
\usepackage{algpseudocodex}
\usepackage{appendix}
\usepackage{amsthm}
\usepackage{amsmath}
\usepackage{makecell}
\usepackage{multirow}
\usepackage{setspace}
\usepackage{subfigure}
\usepackage{natbib}

\hypersetup{colorlinks=true, linkcolor=blue, filecolor=blue, urlcolor=black, citecolor=blue}

\newtheorem{proposition}{Proposition}

\algrenewcommand\algorithmicrequire{\textbf{Input:}}
\algrenewcommand\algorithmicensure{\textbf{Output:}}

\newcommand{\vect}{\mathbf}
\newcommand{\set}{\mathcal}
\newcommand{\dom}{\mathbb}

\newcommand{\indicator}{\mathbb{I}}

\newcommand{\algand}{\textbf{and }}

\title{A Hierarchical Constructive Heuristic for Large-Scale Survivable Traffic Grooming Problem under Double-Link Failures}

\author{
Silong Zhang\\
Department of Logistics and Maritime Studies\\
Hong Kong Polytechnic University\\
\And
Jixuan Feng\\
Department of Logistics and Maritime Studies\\
Hong Kong Polytechnic University\\
\And
Junyan Liu\\
Theory Lab, 2012 Labs\\
Huawei Technologies, Co. Ltd\\
\And
Yu Liu\\
Department of Logistics and Maritime Studies\\
Hong Kong Polytechnic University\\
\And
Zhou Xu\footnotemark\\
Department of Logistics and Maritime Studies\\
Hong Kong Polytechnic University\\
\And
Fan Zhang\\
Theory Lab, 2012 Labs\\
Huawei Technologies, Co. Ltd\\
}

\begin{document}

\renewcommand{\thefootnote}{}
\footnotetext{*Corresponding author: \texttt{lgtzx@polyu.edu.hk}}

\maketitle

\begin{abstract}
This paper studies a survivable traffic grooming problem in large-scale optical transport networks under double-link failures (STG2). Each communication demand must be assigned a route for every possible scenario involving zero, one, or two failed fiber links. Protection against double-link failures is crucial for ensuring reliable telecommunications services while minimizing equipment costs, making it essential for telecommunications companies today. However, this significantly complicates the problem and is rarely addressed in existing studies. Furthermore, current research typically examines networks with fewer than 300 nodes, much smaller than some emerging networks containing thousands of nodes. To address these challenges, we propose a novel hierarchical constructive heuristic for STG2. This heuristic constructs and assigns routes to communication demands across different scenarios by following a hierarchical sequence. It incorporates several innovative optimization techniques and utilizes parallel computing to enhance efficiency. Extensive experiments have been conducted on large-scale STG2 instances provided by our industry partner, encompassing networks with 1,000 to 2,600 nodes. Results demonstrate that within a one-hour time limit and a 16 GB memory limit set by the industry partner, our heuristic improves the objective values of the best-known solutions by 18.5\% on average, highlighting its significant potential for practical applications.
\end{abstract}

\keywords{Traffic grooming \and double-link failure \and heuristic \and optical transport network \and telecommunications}

\section{Introduction}
\label{sect:intro}

\emph{Traffic grooming} optimizes \emph{optical transport networks} (OTNs) by consolidating low-rate traffic flows into high-rate \emph{lightpaths}, reducing facility costs and increasing network throughput \citep{Wu2020}. In this context, each demand must be assigned a route connecting its two terminal nodes, subject to various side constraints. This is crucial given the ever-increasing global telecommunications demands \citep{Cisco2024}. To ensure service reliability, telecommunications companies often require protection against double-link failures, maintaining communication routes for demands even when two fiber links fail. This motivates our study on \emph{survivable traffic grooming under double-link failures} (STG2).

In the existing literature on traffic grooming, double-link failures are rarely considered, and the network sizes examined are of no more than 300 nodes, significantly smaller than some emerging real-world networks involving thousands of nodes. In this paper, we aim to tackle the STG2 in large-scale networks and consider several practical side constraints. These constraints facilitate fast switching between working routes under the scenario without link failures and backup routes under the scenarios with link failures. Solving the STG2 problem is challenging due to the complexity introduced by the large-scale network, numerous failure scenarios, and practical constraints involving interdependency among scenarios. Additionally, practitioners often impose time (e.g., one hour) and memory limits (e.g., 16 GB) on the solution method to ensure efficiency. To address these challenges, we develop a novel hierarchical constructive heuristic for the STG2 in this study.


\subsection{Background of Survivable Traffic Grooming}

The traffic grooming problem is based on an OTN, a telecommunications network that transmits high-speed, high-capacity, and reliable traffic flows using optical fibers. It serves as a critical backbone for internet connections across extensive geographical areas. An OTN consists of interconnected nodes, such as data centers and switching centers, and links made of optical fibers that enable high-speed transmission over long distances. The network manages a series of demands, each defined by a specific traffic volume and two terminal nodes. These demands are strategically routed to minimize facility costs and ensure service reliability.

\begin{figure}[htbp]
	\centering
	\subfigure[Utilizing WDM, the optical signal on a link is divided into 3 wavelengths, capable of transmitting traffic flows independently.]{
		\label{fig:WDM}
		\includegraphics[width = 0.45 \textwidth]{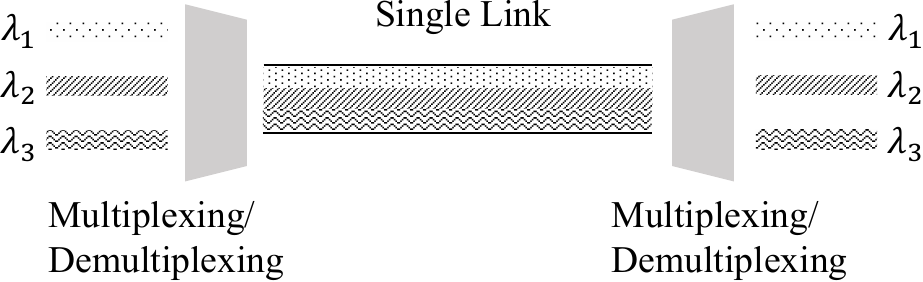}
	}
	\hspace{20pt}
	\subfigure[Each wavelength has a bandwidth capacity of 100 Gbps. Multiple demands can share a wavelength, provided the total traffic does not exceed this capacity.]{
		\label{fig:traffic-grooming}
		\includegraphics[width = 0.45 \textwidth]{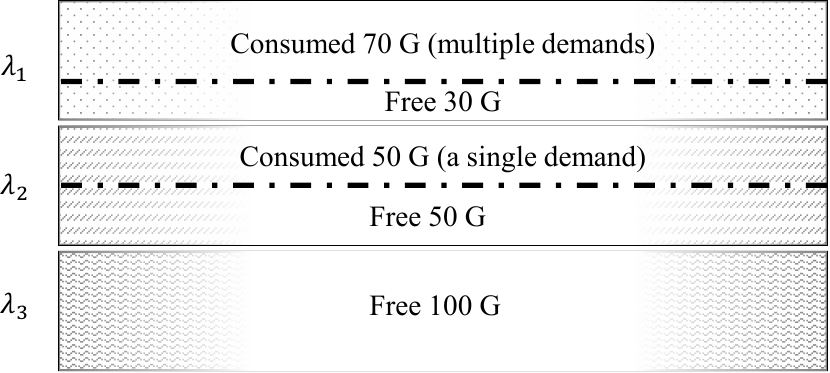}
	}
	\caption{Link capacity in a WDM network enabling traffic grooming. \label{fig:WDM-TG}}
\end{figure}

Traffic grooming enables higher utilization of network resources, compared to the widely studied \emph{routing and wavelength assignment} (RWA) problem. Both problems are based on OTNs utilizing \emph{wavelength division multiplexing} (WDM), which divides the optical signal on a single link into multiple wavelengths (see Figure~\ref{fig:WDM}), each capable of transmitting traffic flows independently. In the RWA, a wavelength on a link is assigned to at most one demand. In the traffic grooming problem, the network throughput can be further increased, as multiple demands can share the same wavelength on a link, provided that their combined traffic does not exceed the bandwidth capacity (see Figure~\ref{fig:traffic-grooming}). As a result, it can be seen that the RWA is a special case of the traffic grooming problem.

Compared to RWA, the traffic grooming problem requires a significantly greater number of decisions, as it must determine how different communication demands share the lightpaths. A lightpath, defined by a wavelength and a sequence of end-to-end connected links, serves as a communication channel between two endpoints equipped with signal processing devices. Multiple low-speed demands can be multiplexed onto this lightpath and transmitted between the endpoints using the designated wavelength and links. It aims to minimize the total number of established lightpaths, which is proportional to the number of required devices, constituting a major portion of network costs \citep{Wu2015}.

\begin{figure}[htbp]
	\centering
	\subfigure[Routing in the working scenario.]{
		\label{fig:scenario-working}
		\includegraphics[width = 0.4 \textwidth]{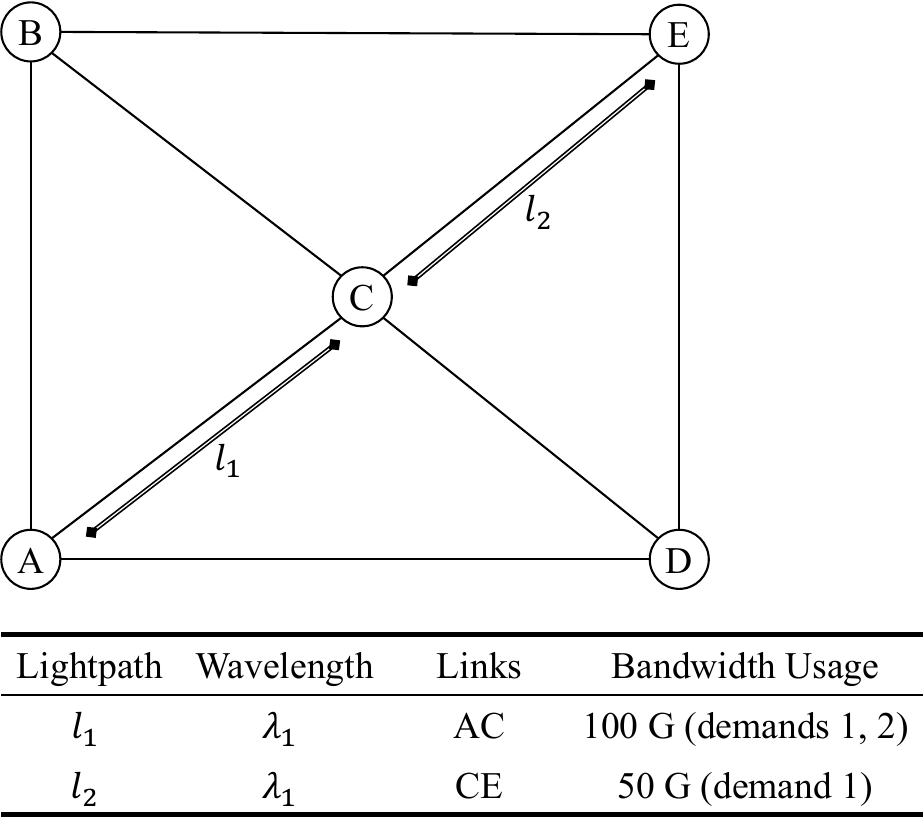}
	}
	\hspace{20pt}
	\subfigure[Routing in the failure scenario where link AC fails.]{
		\label{fig:scenario-AC}
		\includegraphics[width = 0.4 \textwidth]{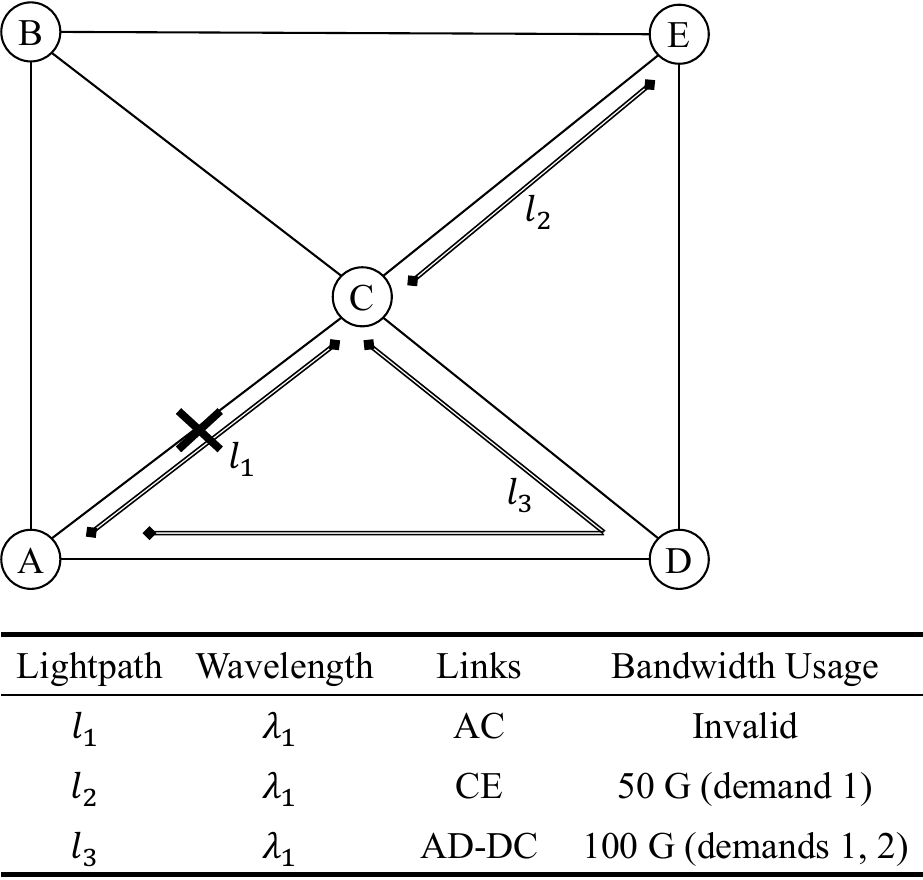}~~~~~
	}
	\\
	\subfigure[Routing in the failure scenario where links AC and CE fail successively.]{
		\label{fig:scenario-AC-CE}
		\includegraphics[width = 0.7 \textwidth]{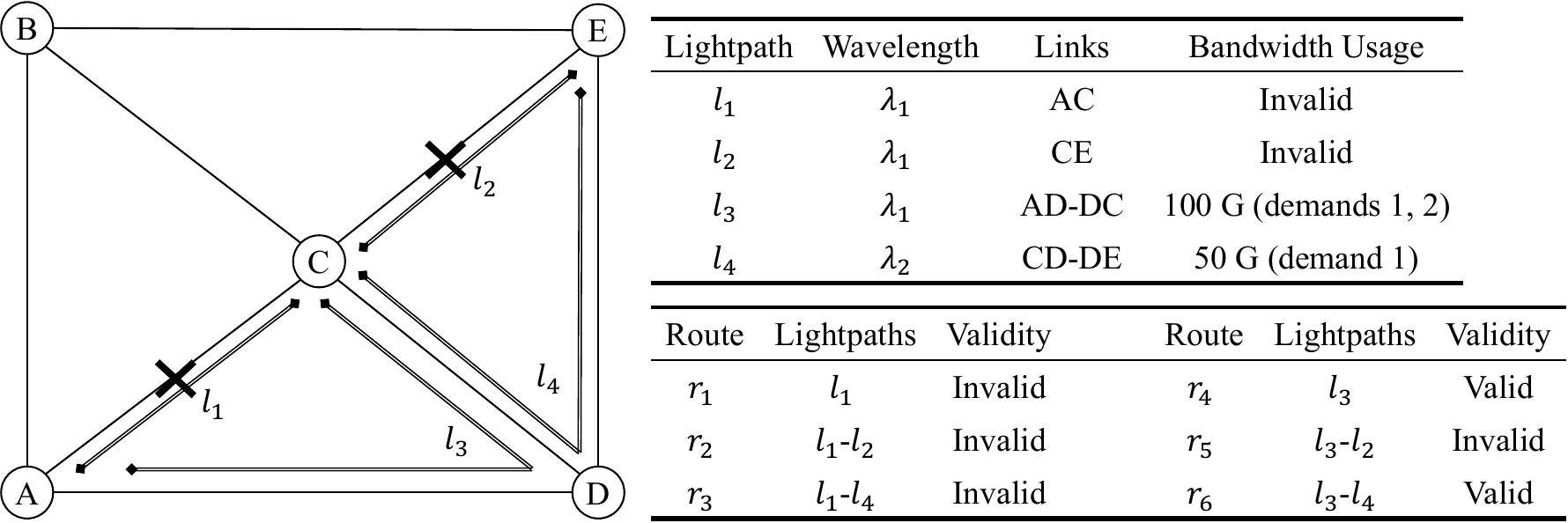}
	}
	\caption{Routing in 3 scenarios for an OTN with 5 nodes, 8 links, and 2 demands: (A, E, 50 G), (A, C, 50 G). \label{fig:intro:scenario}}
\end{figure}

The STG2 enhances network reliability by ensuring survivability against link failures, such as accidental fiber cuts that disrupt routes. A double-link failure occurs when one link fails and another fails before the first is repaired. As OTN sizes continue to expand dramatically, the probability of double-link failures has increased \citep{Sasithong2020}. To improve reliability, STG2 considers the working scenario without link failures and all scenarios involving a single link failure or two link failures, as required by many telecommunications companies today. This involves planning a route for each demand in every scenario, ensuring that backup routes in failure scenarios do not traverse failed links, thereby satisfying demands despite link failures. The number of considered failure scenarios is the square of the number of links. For example, in an OTN with 8 links shown in Figure~\ref{fig:intro:scenario}, there are 65 scenarios: 1 working scenario, 8 single-link failure scenarios, and 56 double-link failure scenarios. Figure~\ref{fig:intro:scenario} illustrates routing plans in 3 scenarios, indicating that lightpaths $l_3$ and $l_4$ should be pre-established for backup routes, although they are not used in the working scenario.

\subsection{Challenges for Solving Large-Scale STG2 and Justifications for Developing Constructive Heuristic}
\label{sec:challenges}

The classic traffic grooming problem is NP-hard \citep{Zhu2002, Wang2008}, and the network size significantly impacts the efficiency of solution methods. In related studies, only instances of no more than 300 nodes are solved (see Table \ref{tab:scale}). We aim to tackle large-scale STG2 with thousands of nodes, aligning with evolving practical requirements driven by the rapid expansion of OTNs. The large network sizes result in more decisions to make, posing challenges for optimizing routes for communication demands.

\begin{table}[b]
\centering
\caption{Largest scales of RWA and traffic grooming instances solved in the literature \label{tab:scale}}
{\begin{tabular}{lllrr}
\hline
Reference & Problem & Failure Type & \# Node & \# Link \\
\hline
\citet{Jaumard2017} & RWA & failure-free & 90 & 350 \\
\citet{Sebbah2012} & RWA & single-link & 43 & 71 \\
\citet{Feng2010} & RWA & double-link & 47 & 98 \\
\citet{Zhu2003} & traffic grooming & failure-free & 277 & unknown \\
\citet{Niu2024} & traffic grooming & single-link & 227 & 409 \\
\citet{Assi2006} & traffic grooming & double-link & 24 & 43 \\
\hline
\end{tabular}}
\end{table}

Another challenge originates from the large number of failure scenarios considered. For survivable OTN optimization problems under single-link (double-link, respectively) failures, the number of failure scenarios is equal to the number of links (the square of the number of links, respectively). This, together with Table \ref{tab:scale}, indicates that the number of failure scenarios considered in the literature is typically no more than 10,000. However, for the large-scale STG2 instances examined in this study, involving thousands of links, there are millions of failure scenarios to consider. This results in significantly more decisions to make, constraints to verify, and data to store, thereby complicating the problem.

Our study also incorporates several practical constraints that challenge algorithm design. For instance, consistent routing and bandwidth reservation constraints are essential for fast switching between working and backup routes, ensuring low-latency connections. These constraints create interdependency in routing and bandwidth consumption across scenarios, significantly complicating the algorithm design.

Limited computing resources, such as computing time and available memory, further heighten the challenge of algorithm design. Achieving short computing time is essential for telecommunications service providers to optimize subsequent decisions based on the STG2 solutions. Additionally, major providers often run multiple computing tasks simultaneously on their platforms, restricting the memory and time available for solving STG2.

In this study, a one-hour time limit and a 16 GB memory limit are set by our industry partner, Huawei Technologies. Preliminary experiments indicate that for their large-scale STG2 instances with thousands of nodes, the one-hour time limit is nearly reached just by verifying a solution's feasibility. Additionally, importing linear programming relaxation models of these instances far exceeds the 16 GB memory limit. Consequently, methods based on iterative search or mathematical programming are not applicable. Therefore, we focus on developing a constructive heuristic.

\subsection{Related Works}
\label{sec:review}

\subsubsection*{Traffic Grooming without Considering Failures}
\label{sect:review:TG}

Constructive heuristics have been developed for the traffic grooming problem without considering failures. \citet{Zhu2003} and \citet{Yao2005a} construct solutions by sequentially planning routes for each demand. \citet{Zhu2003} introduce an auxiliary graph with a logical layer for lightpaths and multiple wavelength layers for wavelength-link pairs. \citet{Yao2005a} propose a link-bundled auxiliary graph (LBAG), which consolidates all wavelength layers into a single physical layer. Both studies use shortest-path algorithms to construct routes. However, these are not applicable to our STG2 due to various side constraints that must be incorporated.

\citet{Chen2010} employ a three-phase method to construct solutions and groom traffic within a hub-and-spoke framework. They begin by decomposing the network into clusters, then determine endpoints of lightpaths for grooming intracluster and intercluster traffic, and finally solve an RWA problem to decide the physical routing and wavelength of lightpaths. When adapted to our problem, the third phase may become infeasible due to various side constraints and interdependency between scenarios. \citet{Wu2015} construct a solution by sequentially finding routes for each demand, followed by optimizing the logical and physical routing separately to enhance the solution. However, they do not address optical reach constraints or limit the number of wavelengths, which are crucial in practice.

Heuristics based on mathematical programming and iterative search have also been applied to traffic grooming problems. \citet{Dawande2007a} use heuristics based on column generation and Lagrangian relaxation, requiring that demands with different terminal nodes do not share the same lightpath. \citet{Vignac2016} apply Benders decomposition and Dantzig–Wolfe reformulation under strict assumptions regarding demand granularity, routing schemes, and grooming configurations. \citet{Wu2020} first construct an initial feasible solution and then iteratively search neighboring solutions to reduce the number of lightpaths. However, applying these methods to large-scale STG2 is too time-consuming and memory-intensive due to the extensive network sizes and numerous failure scenarios.


\subsubsection*{Survivable Traffic Grooming}
\label{sect:review:survivable}

Several different survivability schemes have been proposed for traffic grooming problems considering failures, broadly classified into two categories: Pre-planned protection and dynamic restoration \citep{Somani2008}. Pre-planned protection schemes reserve resources that remain idle during the working scenario, ensuring recovery when failures occur. In contrast, dynamic restoration schemes do not reserve resources beforehand, instead searching for spare resources after a failure. \citet{Assi2006} explore a hybrid scheme for STG2, incorporating protection for single-link failures and restoration for double-link failures, although without providing guaranteed protection for the latter. This study focuses on a protection scheme to enhance network reliability.

According to the literature, various protection schemes can be applied at the lightpath level (PAL) and the connection level (PAC). In PAL (PAC, respectively), all traffic on a lightpath (route, respectively) switches to a backup lightpath (route, respectively) upon failure. Since a route is a concatenation of lightpaths, PAL is a specialization of PAC and cannot achieve higher resource utilization than PAC. To minimize facility costs, we focus on PAC in this study. \citet{Yao2005b} and \citet{Jaekel2008} explore survivable traffic grooming under single-link failures (STG1) using PAL, where each lightpath is accompanied by a link-disjoint backup lightpath. They identify these lightpaths using a shortest-path algorithm and an integer linear programming (ILP) formulation, respectively.

PAC can be further classified into full path protection (FPP) and partial path protection (PPP). FPP requires that working and backup routes are link-disjoint, while PPP does not. \citet{Yao2005c} and \citet{Niu2024} investigate STG1 using FPP. \citet{Yao2005c} focus on maximizing network throughput and propose two grooming heuristics to establish lightpaths and route traffic, either separately or jointly. \citet{Niu2024} aim to minimize the number of lightpaths while ensuring all demands are protected, proposing a heuristic that plans backup routes after establishing working routes for all demands. FPP has also been applied to STG2 \citep{Yang2009} and RWA under double-link failures \citep{Feng2010}. \citet{Jaekel2012} explore STG1 with PPP, presenting ILP formulations and a heuristic method that sequentially plans each failure scenario.

To the best of our knowledge, STG2 with PPP has not been explored in the literature. As a less restrictive approach, PPP serves as a generalization of FPP, offering the potential for higher resource utilization and cost savings. Therefore, we focus on PPP in this study, albeit with higher computational complexity.


\subsection{Contributions}
Existing studies reviewed in Section~\ref{sec:review} typically adopt mathematical programming and iterative search methods. However, they are not applicable to our study due to the large network sizes, extensive constraints, and imposed computing time and memory limits (1 hour and 16 GB).
Therefore, as explained in Section~\ref{sec:challenges}, to effectively address large-scale STG2, we develope an efficient constructive heuristic with the following key contributions:
\begin{itemize}
	\item We propose constructing routes and assigning them to communication demands by following a hierarchical sequence of scenarios. This approach mitigates interdependency between scenarios, enabling our heuristic to efficiently manage consistent routing constraints and validate bandwidth capacity constraints.
	\item We adopt a reuse-then-search routing strategy to efficiently satisfy demands in a large-scale network, considering practical constraints and minimizing objective value increase. It first attempts to reuse existing routes to satisfy new demands without objective value increase. If it fails, a tailored labeling algorithm is employed to search for a new route.
	\item To manage the large number of demands and simplify decision-making, we aggregate multiple demands into a single one, assigning them the same routing plan. To minimize the number of aggregated demands while respecting capacity constraints, we introduce a bin-packing problem for each pair of terminal nodes and solve it with an efficient heuristic.
	\item To efficiently manage the large number of double-link failure scenarios, we develop a parallel algorithm based on our constructive heuristic, utilizing multiple threads to concurrently plan distinct scenarios. We implement mechanisms to prevent the generation of substitutable lightpaths, mitigate race conditions (where multiple threads concurrently access and modify shared data), and minimize thread idleness. To the best of our knowledge, this is the first parallel algorithm for survivable OTN problems.
	\item We conduct extensive computational experiments on large-scale STG2 instances provided by our industrial partner. The results demonstrate the effectiveness of our constructive heuristic and various enhancements. Within a one-hour time limit and a 16 GB memory limit, our heuristic can solve instances with over 2,000 nodes, 3,000 links, 10,000 demands, and 9 million scenarios, improving the objective value of the best-known solutions by 18.5\% on average, highlighting its significant potential for practical applications.
\end{itemize}

We next present the problem statement in Section~\ref{sect:problem}, the hierarchical constructive heuristic in Section~\ref{sect:construct}, and the enhancements by parallel computing in Section~\ref{sect:parallel}. Our computational results are reported in Section~\ref{sect:experiments}, followed by the conclusion in Section~\ref{sect:conclusion}.

\section{Problem Description of STG2}
\label{sect:problem}

The STG2 is defined on an OTN represented by an undirected graph $(\set{N}, \set{E})$, where $\set{N}$ is the set of nodes and $\set{E}$ is the set of undirected links. Each link $e \in \set{E}$ has an associated length $d(e)$ and a set of wavelengths $\Lambda$. Each wavelength $\lambda \in \Lambda$ has a bandwidth capacity $B$. There is a set of demands $\set{K}$, where each demand $k \in \set{K}$ has a specific traffic volume $b_k \le B$, requiring transmission between two terminal nodes $s_k \in \set{N}$ and $t_k \in \set{N}$. Figure~\ref{fig:intro:scenario} illustrates an example featuring 5 nodes, 8 links, 2 wavelengths, and 2 demands.

The STG2 involves route planning for demands across all possible scenarios involving zero, one, or two failed links. The working scenario (no link failures) is denoted as $(0, 0)$, the \emph{level-1 scenario} where a single link $e_1$ fails is denoted as $(e_1, 0)$, and the \emph{level-2 scenario} where two links $e_1$ and $e_2$ fail successively is denoted as $(e_1, e_2)$. The set of all considered scenarios is represented by $\Omega = \{(0, 0)\} \cup \Omega^{1} \cup \Omega^{2}$, where $\Omega^{1} = \{(e_1, 0): e_1 \in \set{E}\}$ and $\Omega^{2} = \{(e_1, e_2): e_1, e_2 \in \set{E}, e_1 \ne e_2\}$.

The STG2 aims to minimize the number of established lightpaths, ensuring that each demand is assigned a route in every considered scenario, without traversing failed links. A lightpath $l = (\lambda, (e_{l}^{1}, e_{l}^{2}, \ldots, e_{l}^{n}))$, defined by a wavelength $\lambda \in \Lambda$ and a sequence of distinct, end-to-end connected links $e_{l}^{1}, e_{l}^{2}, \ldots, e_{l}^{n}$, serves as a communication channel between two endpoints. This definition implies that a lightpath $l$ has a bandwidth capacity of $B$ and must occupy the same wavelength $\lambda$ on traversed links $e_{l}^{1}, e_{l}^{2}, \ldots, e_{l}^{n}$, known as the \emph{wavelength continuity constraints} \citep{Kennington2003}. For example, the lightpath $l_3 = (\lambda_1, (\text{AD, CD}))$ shown in Figure~\ref{fig:intro:scenario} occupies the same wavelength $\lambda_1$ on links AD and CD, and can transmit up to 100 G of traffic between endpoints A and C. Additionally, due to signal attenuation, the \emph{optical reach constraints} \citep{Chen2015, Yildiz2017} require that the length of a lightpath $l$ does not exceed the optical reach limit $\bar{d}$, i.e., $\sum_{e \in \set{E}(l)} d(e) \le \bar{d}$ where $\set{E}(l) = \{e_{l}^{1}, e_{l}^{2}, \ldots, e_{l}^{n}\}$ is the set of links traversed by $l$. A route $r = (l_{r}^{1}, l_{r}^{2}, \ldots, l_{r}^{n})$ is a sequence of distinct, end-to-end connected lightpaths. For example, $r_2 = (l_1, l_2)$, depicted in Figure~\ref{fig:intro:scenario}, is a route connecting nodes A and E, capable of transmitting demand $k_1$. Denote $\set{L}$ and $\set{R}$ as the set of all lightpaths and routes, respectively. Denote $r_{k}^{e_1 e_2}$ as the route for demand $k \in \set{K}$ used in scenario $(e_1, e_2) \in \Omega$.

The STG2 also considers some other side constraints. \emph{Simple path constraints} require a route to be elementary, meaning it traverses each node at most once. For example, in Figure~\ref{fig:intro:scenario}, routes $r_1, r_2, r_3, r_4, r_5$ are elementary, whereas $r_6$ is not, as it traverses node D twice. Although these constraints can reduce system disorders \citep{Wu2015}, they may increase the objective value. Therefore, we require only that working routes ($r_{k}^{0 0}, k \in \set{K}$) be elementary. In a failure scenario, each demand should be transmitted along a route $r$ without traversing any failed link, i.e., $e_1 \not\in \set{E}(r)$ in scenario $(e_1, 0) \in \Omega^{1}$ and $e_1, e_2 \not\in \set{E}(r)$ in scenario $(e_1, e_2) \in \Omega^{2}$, where $\set{E}(r)$ is the set of links traversed by route $r$. \emph{Wavelength assignment constraints} are widely considered for OTNs \citep{Daryalal2022}, requiring that a wavelength on a link be occupied by at most one established lightpath. In other words, any two lightpaths $l$ and $l'$ using the same wavelength ($\lambda(l) = \lambda(l')$, where $\lambda(l)$ represents the wavelength occupied by lightpath $l$) cannot traverse a common link ($\set{E}(l) \cap \set{E}(l') = \emptyset$). For example, in Figure~\ref{fig:intro:scenario}, lightpaths $l_3$ and $l_4$ must occupy different wavelengths ($\lambda_1$ and $\lambda_2$), because they share the link CD. Denote $\set{L}(r)$ as the set of lightpaths included in route $r$. Note that different lightpaths in $\set{L}(r)$ can occupy different wavelengths.

\begin{figure}[t]
	\centering
	\includegraphics[width = 0.6 \textwidth]{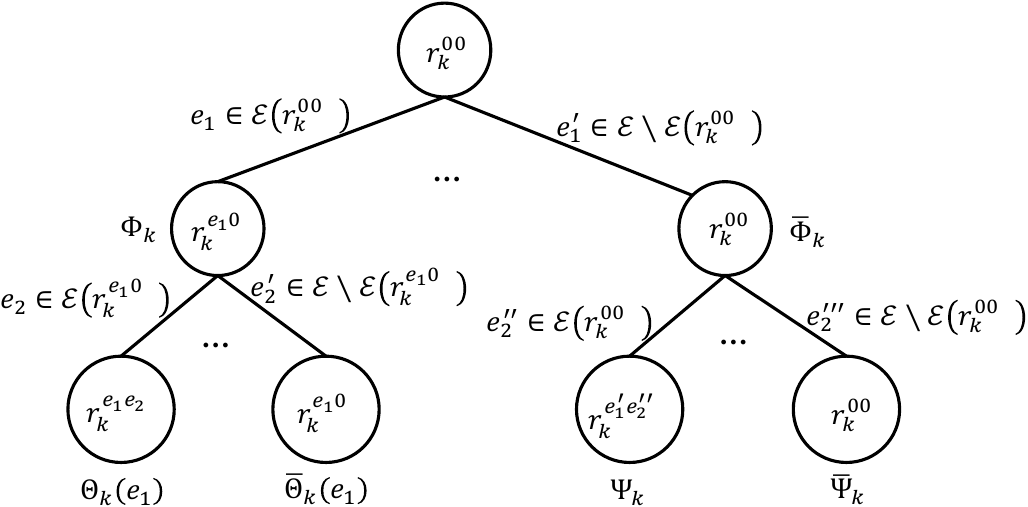}
	\caption{Route planning for demand $k \in \set{K}$ in scenarios $\Omega$, where failed links are represented on the edges. \label{fig:problem:route-planning}}
\end{figure}

\emph{Consistent routing constraints} are essential for fast route switching and low-latency communications, as they minimize the number of routes that need to be altered when a link fails. These constraints ensure that the route used before a failure continues to be used if it does not traverse the failed link. Specifically, for each demand $k \in \set{K}$, the working route $r_{k}^{0 0}$ should remain in use in level-1 scenarios $\bar{\Phi}_{k} = \{(e_1, 0) \in \Omega^{1}: e_1 \in \set{E} \setminus \set{E}(r_{k}^{0 0})\}$ and level-2 scenarios $\bar{\Psi}_{k} = \{(e_1, e_2) \in \Omega^{2}: e_1, e_2 \in \set{E} \setminus \set{E}(r_{k}^{0 0})\}$. Similarly, for each link $e_1 \in \set{E}(r_{k}^{0 0})$, the level-1 route $r_{k}^{e_1 0}$ should remain in use in level-2 scenarios $\bar{\Theta}_{k}(e_1) = \{(e_1, e_2) \in \Omega^{2}: e_2 \in \set{E} \setminus \set{E}(r_{k}^{e_1 0})\}$. Consequently, the routing plan for a demand $k \in \set{K}$ forms a tree structure, as illustrated in Figure~\ref{fig:problem:route-planning}.

\emph{Bandwidth reservation constraints} are also essential for fast route switching and low-latency communications. When failed links are repaired, all demands revert to their working routes, ensuring protection against future double-link failures. These constraints require that the bandwidth on lightpaths consumed by demand $k \in \set{K}$ in the working scenario be reserved. The reserved bandwidth can be reused by demand $k$ in failure scenarios, but cannot be used by any other demand in any scenario. This ensures that demand $k$ can revert to the working route without waiting for other demands to release bandwidth resources, thereby accelerating route switching. Consequently, the \emph{bandwidth capacity constraints} require that the total bandwidth consumption, including reserved bandwidth, on each lightpath in any scenario must not exceed $B$.

The STG2 can be formulated as an ILP (see Appendix~\ref{sect:formulation}). However, we do not utilize it to develop solution methods, because for large-scale STG2 instances with thousands of nodes examined in this study, importing their linear programming relaxation models far exceeds the 16 GB memory limit set by our industry partner.

\section{Hierarchical Constructive Heuristic}
\label{sect:construct}

To efficiently solve large-scale STG2 instances (thousands of nodes and millions of failure scenarios) with limited computing time and memory, we propose a novel hierarchical constructive heuristic, outlined in Section~\ref{sect:construct:outline}. Section~\ref{sect:demand:agg} illustrates the method for aggregating demands. In Section~\ref{sect:labeling}, we develop a labeling algorithm for finding routes, and Section~\ref{sect:reuse-route} introduces a technique for reusing generated routes. Our  hierarchical constructive heuristic is further enhanced by efficient checks of bandwidth capacity constraints (see Appendix~\ref{sect:bandwidth:check}) and efficient updates of bandwidth consumption (see Appendix~\ref{sect:bandwidth:update}).

\subsection{Algorithm Outline}
\label{sect:construct:outline}

The hierarchical constructive heuristic is outlined in Algorithm~\ref{alg:construction:outline}. It begins by aggregating demands (line~\ref{alg:construction:outline:agg:demand}) to reduce the problem size (Section~\ref{sect:demand:agg}), followed by sorting the sequence of demands (line~\ref{alg:construction:outline:sort:demand}). The solution is constructed hierarchically: First plan the working scenario (lines~\ref{alg:construction:outline:working:begin}--\ref{alg:construction:outline:working:end}), then plan level-1 scenarios (lines~\ref{alg:construction:outline:level-1:begin}--\ref{alg:construction:outline:level-1:end}), and finally plan level-2 scenarios (lines~\ref{alg:construction:outline:level-2:begin}--\ref{alg:construction:outline:level-2:end}). In each scenario, routes for demands are planned sequentially. Working routes are determined using the labeling algorithm (Section~\ref{sect:labeling}).  When planning level-1 and level-2 routes, to reduce the number of calls of the labeling algorithm, we first try to reuse routes generated during the solution process (see Section~\ref{sect:reuse-route}). If this fails, we then employ the labeling algorithm to generate a new route.

\begin{algorithm}
\begin{small}
\caption{Hierarchical Constructive Heuristic}
\label{alg:construction:outline}
\begin{algorithmic}[1]
	\Statex \algorithmicrequire \text{ } $\set{N}, \set{E}, \Lambda, \set{K}, B, \bar{d}$
	\Statex \algorithmicensure \text{ } Routes $r_{k}^{e_1 e_2}, k \in \set{K}, (e_1, e_2) \in \Omega$
	\State Aggregate demands $\set{K}$ into $\set{K}'$ \label{alg:construction:outline:agg:demand} $\qquad \qquad \qquad \qquad \qquad \qquad \qquad \qquad \qquad \qquad \qquad \qquad$ (Section~\ref{sect:demand:agg})
	\State Sort demands $\set{K}'$ in non-increasing order of traffic volume \label{alg:construction:outline:sort:demand}
	\ForAll{demand $k \in \set{K}'$} \label{alg:construction:outline:working:begin}
		\State Find the working route $r_{k}^{0 0}$ \label{alg:construction:outline:working:labeling} with the labeling algorithm $\qquad \qquad \qquad \qquad \qquad \quad$ (Section~\ref{sect:labeling})
		\State $\Phi_{k} \gets \{(e_1, 0): e_1 \in \set{E}(r_{k}^{0 0})\}$; $\Psi_{k} \gets \{(e_1, e_2): e_1 \in \set{E} \setminus \set{E}(r_{k}^{0 0}), e_2 \in \set{E}(r_{k}^{0 0})\}$ \label{alg:construction:outline:working:end}
	\EndFor
	\State Sort scenarios $\Omega^{1}$ \label{alg:construction:outline:sort:scn:1} in non-increasing order of the number of unplanned demands
	\ForAll{scenario $(e_1, 0) \in \Omega^{1}$} \label{alg:construction:outline:level-1:begin}
		\ForAll{demand $k \in \set{K}'$ where $(e_1, 0) \in \Phi_{k}$}
			\State Find the level-1 route $r_{k}^{e_1 0}$ to protect $k$ in scenario $(e_1, 0)$ by reusing an existing route (Section~\ref{sect:reuse-route}). If this fails, use the labeling algorithm (Section~\ref{sect:labeling}) to find a new route \label{alg:construction:outline:level-1:route}
			\State $\Theta_{k}(e_1) \gets \{(e_1, e_2): e_2 \in \set{E}(r_{k}^{e_1 0})\}$ \label{alg:construction:outline:level-1:end}
		\EndFor
	\EndFor
	\State Sort scenarios $\Omega^{2}$ in non-increasing order of the number of unplanned demands \label{alg:construction:outline:sort:scn:2}
	\ForAll{scenario $(e_1, e_2) \in \Omega^{2}$} \label{alg:construction:outline:level-2:begin}
		\ForAll{demand $k \in \set{K}'$ where $(e_1, e_2) \in \Psi_{k} \cup \Theta_{k}(e_{1})$}
			\State Find the level-2 route $r_{k}^{e_1 e_2}$ to protect $k$ in scenario $(e_1, e_2)$ by reusing an existing route (Section~\ref{sect:reuse-route}). If this fails, use the labeling algorithm (Section~\ref{sect:labeling}) to find a new route \label{alg:construction:outline:level-2:route} \label{alg:construction:outline:level-2:end}
		\EndFor
	\EndFor
\end{algorithmic}
\end{small}
\end{algorithm}

The hierarchical constructive heuristic offers several advantages. Firstly, it naturally satisfies consistent routing constraints by sequentially planning the working route, level-1 routes, and level-2 routes. Specifically, once demand $k \in \set{K}'$ is assigned with a working route $r_{k}^{0 0}$, this route is also assigned to the demand in level-1 scenarios $\bar{\Phi}_{k}$ and level-2 scenarios $\bar{\Psi}_{k}$. Similarly, when demand $k \in \set{K}'$ is assigned with a level-1 route $r_{k}^{e_1 0}$ in scenario $(e_1, 0) \in \Phi_{k}$, this route is also assigned in level-2 scenarios $\bar{\Theta}_{k}(e_1)$ by reserving bandwidth $b_k$ on lightpaths $\set{L}(r_{k}^{e_1 0})$ in scenarios $\bar{\Theta}_{k}(e_1)$. Thus, we only need to explicitly plan the working scenario, level-1 scenarios $\bigcup_{k \in \set{K}'} \Phi_{k}$, and level-2 scenarios $\bigcup_{k \in \set{K}'} \Psi_{k} \cup \Theta_{k}$, where $\Theta_{k} = \bigcup_{e_1 \in \set{E}(r_{k}^{0 0})} \Theta_{k}(e_1)$.

Secondly, while planning a route in one scenario can affect bandwidth consumption in others (due to bandwidth reservation and consistent routing constraints), the hierarchical constructive heuristic allows us to verify bandwidth capacity constraints across all scenarios by checking only the scenario being currently planned (see Appendix~\ref{sect:bandwidth:check}). It also enables updating bandwidth consumption with fewer operations and less memory than a direct update approach (see Appendix~\ref{sect:bandwidth:update}).

Moreover, once the working and all level-1 scenarios are planned, planning a route in a level-2 scenario imposes no constraints on other level-2 scenarios, allowing us to plan these scenarios in parallel (see Section~\ref{sect:parallel}).

\subsection{Demand Aggregation}
\label{sect:demand:agg}

To reduce the problem size, we aggregate demands so that each batch is assigned the same routing plan. Demands $k_1, k_2, \ldots, k_n$ can be aggregated into a single demand $k'$ if and only if they share the same terminal nodes, i.e., $\{s_{k_1}, t_{k_1}\} = \{s_{k_2}, t_{k_2}\} = \cdots = \{s_{k_n}, t_{k_n}\}$, and their total traffic volume does not exceed the bandwidth capacity of a wavelength, i.e., $b_{k_1} + b_{k_2} + \cdots + b_{k_n} \le B$. For the aggregated demand $k'$, we have $s_{k'} = s_{k_1}, t_{k'} = t_{k_1}, b_{k'} = b_{k_1} + b_{k_2} + \cdots + b_{k_n}$, and demands $k_1, k_2, \ldots, k_n$ use the same route as $k'$ in each scenario, i.e., $r_{k_1}^{e_1 e_2} = r_{k_2}^{e_1 e_2} = \cdots = r_{k_n}^{e_1 e_2} = r_{k'}^{e_1 e_2}, \forall (e_1, e_2) \in \Omega$.

For each pair of terminal nodes $\{s, t\}$, decisions on aggregating demands $\set{K}_{s t} = \{k \in \set{K}: \{s_k, t_k\} = \{s, t\}\}$ are made by solving a bin-packing problem. Each demand $k \in \set{K}_{s t}$ has a weight $b_k$, and each bin has a capacity of $B$. The goal is to pack demands $\set{K}_{s t}$ using the minimum number of bins, thereby reducing the number of aggregated demands and the number of routing decisions. To efficiently solve this bin-packing problem, we first sort demands $\set{K}_{s t}$ in non-increasing order of weights, then apply the well-known first-fit policy
: For each demand, place it in the first bin with sufficient remaining capacity; If no such bin exists, place it in a new bin. Demands within a bin are aggregated into a single demand, and we denote $\set{K}'$ as the set of demands after aggregation.

\subsection{Labeling Algorithm for Route Planning}
\label{sect:labeling}


For a demand $k \in \set{K}'$ and scenario $(e_1, e_2) \in \Omega$, the labeling algorithm seeks to search a route $r_{k}^{e_1 e_2}$ that can transmit demand $k$ in scenario $(e_1, e_2)$ while satisfying side constraints. Denote $\set{L}'$ as the set of established lightpaths, and $y_{k l}^{e_1 e_2}, k \in \set{K}', l \in \set{L}', (e_1, e_2) \in \Omega: e_1, e_2 \not\in \set{E}(l)$ as the binary variable that equals one if and only if lightpath $l$ is used by demand $k$ in scenario $(e_1, e_2)$, i.e., $l \in \set{L}(r_{k}^{e_1 e_2})$. Then we know that $C_{l}^{e_1 e_2} = \sum_{k \in \set{K}'} b_k \max\{y_{k l}^{0 0}, y_{k l}^{e_1 e_2}\}$ is the bandwidth consumption on lightpath $l \in \set{L}'$ in scenario $(e_1, e_2) \in \Omega$, according to bandwidth capacity constraints \eqref{eqn:problem:capacity}. At the start of the constructive heuristic, $\set{L}'$ is empty, and both $\vect{y} = (y_{k l}^{e_1 e_2})_{k \in \set{K}', l \in \set{L}', (e_1, e_2) \in \Omega: e_1, e_2 \not\in \set{E}(l)}$ and $\vect{C} = (C_{l}^{e_1 e_2})_{l \in \set{L}', (e_1, e_2) \in \Omega}$ are zero vectors. The underlying graph of the labeling algorithm is detailed in Section~\ref{sect:labeling:topology}, with arc costs defined in Section~\ref{sect:labeling:cost}, label definition and extension developed in Section~\ref{sect:labeling:label}, and heuristic dominance rules provided in Section~\ref{sect:labeling:dominance}.
\begin{equation}
\label{eqn:problem:capacity}
	\sum_{k \in \set{K}'} b_k \max\{y_{k l}^{0 0}, y_{k l}^{e_1 e_2}\} \le B, \quad \forall l \in \set{L}', (e_1, e_2) \in \Omega: e_1, e_2 \not\in \set{E}(l).
\end{equation}

\subsubsection{Two-Layer Graph (LBAG)}
\label{sect:labeling:topology}
To represent the relationships between the physical layer (links and wavelengths) and the logical layer (lightpaths) for STG2, we utilize a two-layer graph, named as LBAG and denoted by $\set{G} = (\set{N} \cup \set{N}', \set{A}_p \cup \set{A}_l \cup \set{A}_t \cup \set{A}_r)$, which was initially proposed by \citet{Yao2005a} for the problem without considering link failures. The physical layer consists of the set of physical nodes $\set{N} = \{1, 2, \ldots, |\set{N}|\}$ and the set of directed physical arcs $\set{A}_p = \{(i_1, i_2, e), (i_2, i_1, e): e \in \set{E}, i_1, i_2 \text{ are endpoints of link } e\}$, where each link $e \in \set{E}$ is associated with two physical arcs. Note that there may be multiple links between a pair of nodes. For each physical arc $a \in \set{A}_p$, denote $e(a)$ as the associated link. The logical layer consists of the set of logical nodes $\set{N}' = \{1', 2', \ldots, |\set{N}|'\}$, where node $i' \in \set{N}'$ is a copy of node $i \in \set{N}$, and the set of directed logical arcs $\set{A}_l = \{(i_{1}^{\prime}, i_{2}^{\prime}, l), (i_{2}^{\prime}, i_{1}^{\prime}, l): l \in \set{L}', i_1, i_2 \text{ are endpoints of lightpath } l\}$, where each lightpath $l \in \set{L}'$ is associated with two logical arcs. Note that there may be multiple lightpaths between a pair of nodes. For each logical arc $a \in \set{A}_l$, denote $l(a)$ as the associated lightpath and $\set{E}(a)$ as the set of links traversed by lightpath $l(a)$. The sets $\set{A}_t = \{(i', i): i \in \set{N}\}$ and $\set{A}_r = \{(i, i'): i \in \set{N}\}$ are transmitter arcs and receiver arcs, respectively, connecting the physical and logical layers.

\begin{figure}[t]
	\centering
	\subfigure[Before routing.]{
		\label{fig:labeling:LBAG:before-routing}
		\includegraphics[width = 0.31 \textwidth]{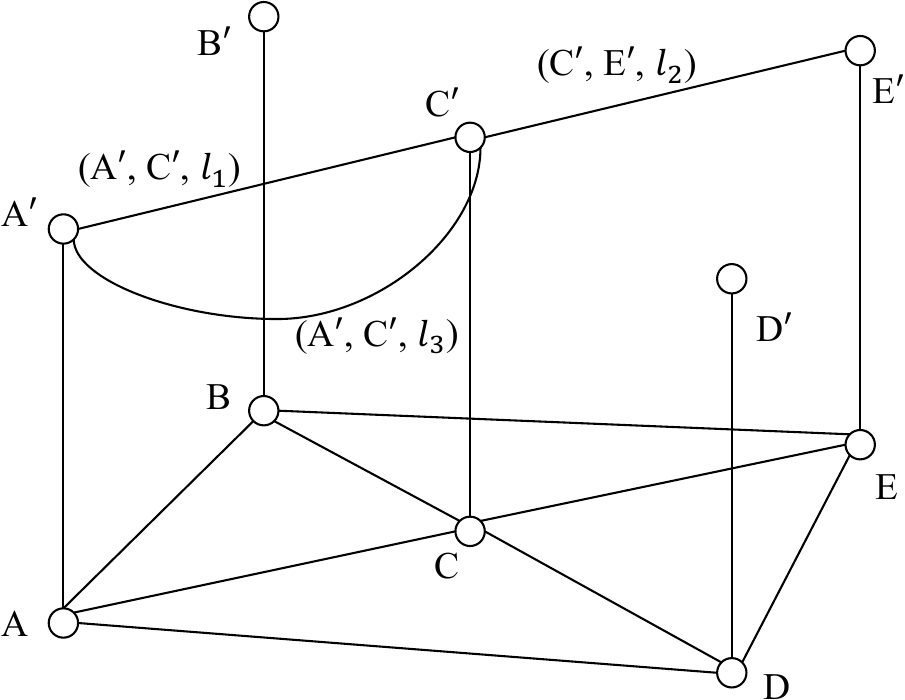}
	}
	\subfigure[A path is found (dotted).]{
		\label{fig:labeling:LBAG:found-route}
		\includegraphics[width = 0.31 \textwidth]{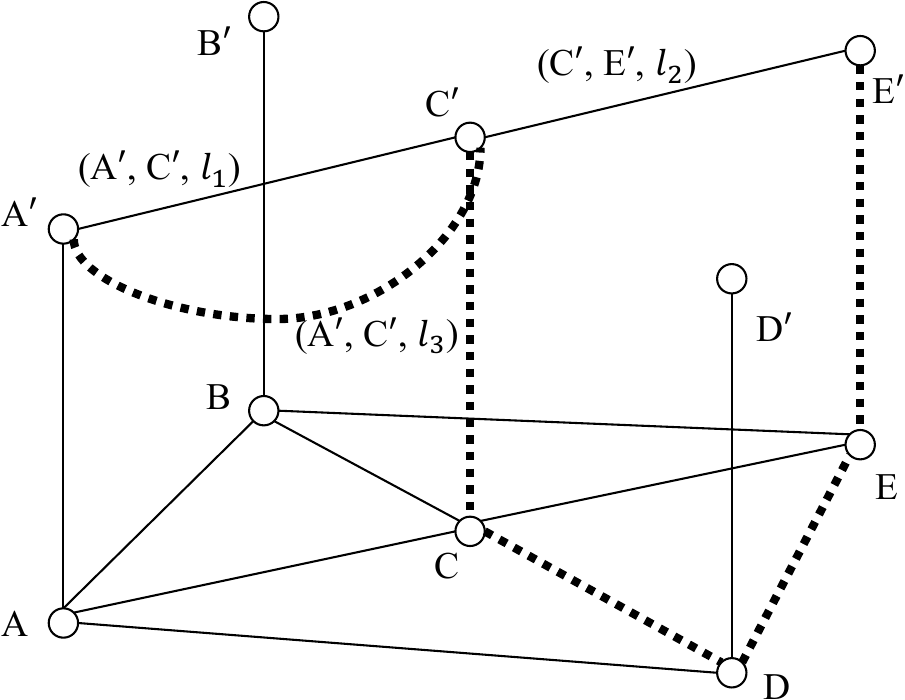}
	}
	\subfigure[Establish a lightpath $l_4$.]{
		\label{fig:labeling:LBAG:establish-lightpath}
		\includegraphics[width = 0.31 \textwidth]{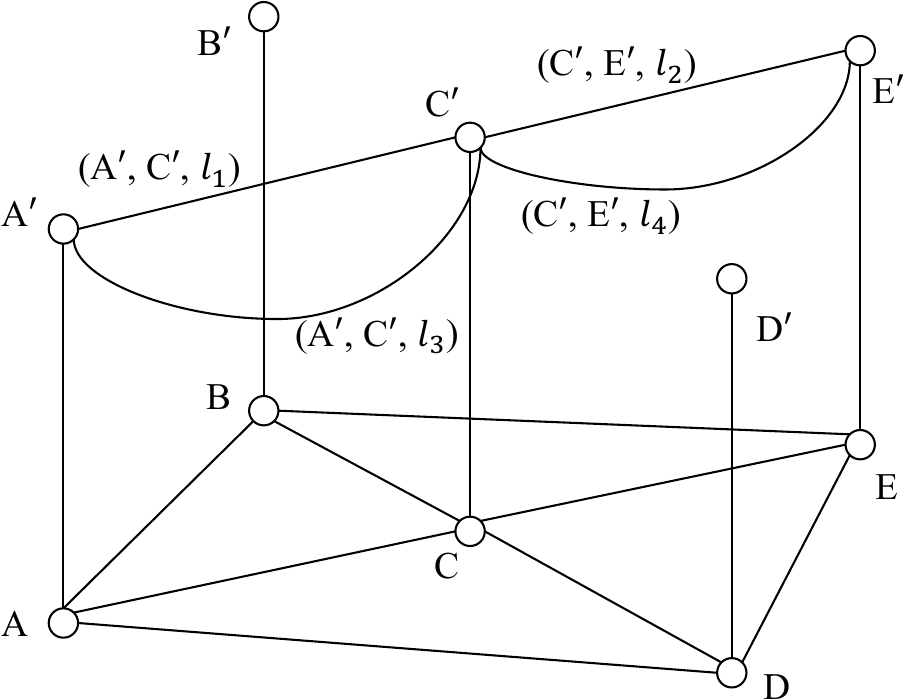}
	}
	\caption{Routing for demand (A, E, 50 G) in scenario (AC, CE), where links AC and CE fail successively. Fugure~\ref{fig:labeling:LBAG:before-routing} shows an LBAG where lightpaths $l_1$ and $l_2$ fail. Fugure~\ref{fig:labeling:LBAG:found-route} shows that the labeling algorithm finds a path sequentially traversing arcs (A$'$, C$'$, $l_3$), (C$'$, C), (C, D), (D, E), (E, E$'$). Fugure~\ref{fig:labeling:LBAG:establish-lightpath} shows that a new lightpath $l_4$ is established, and the demand will use route $(l_3, l_4)$ in the scenario. \label{fig:labeling:LBAG}}
\end{figure}

To find a route for demand $k \in \set{K}'$, the labeling algorithm seeks a path on the LBAG with endpoints $s_{k}^{\prime}$ and $t_{k}^{\prime}$. Traversing a logical arc $a \in \set{A}_l$ indicates reusing an established lightpath $l(a)$. Sequential traversal of a transmitter arc $(i_{1}^{\prime}, i_1) \in \set{A}_t$, physical arcs $a_1, a_2, \ldots, a_n \in \set{A}_p$, and a receiver arc $(i_{n + 1}, i_{n + 1}^{\prime}) \in \set{A}_r$ implies that a new lightpath $(\lambda, (e(a_1), e(a_2), \ldots, e(a_n)))$ should be established, where $\lambda$ is a wavelength unoccupied on these links, i.e., $\lambda \in \bigcap_{m = 1, 2, \ldots, n} \Lambda(e(a_m))$. Here, $\Lambda(e) = \Lambda \setminus \{\lambda(l): l \in \set{L'}, e \in \set{E}(l)\}$ represents the set of unoccupied wavelengths on link $e$. For example, in Figure~\ref{fig:labeling:LBAG}, sequentially traversing arcs (C$'$, C), (C, D), (D, E), (E, E$'$) indicates that a lightpath $l_4$ should be established.

\subsubsection{Arc Costs and Lower Bounds on Path Costs}
\label{sect:labeling:cost}
To effectively guide the labeling algorithm to search for routes of good quality, we need to compute an estimated cost for each label, which will be explained in Section~\ref{sect:labeling:label}. To achieve this, we need to set arc costs according to Table~\ref{tab:labeling:cost}, where $\alpha, \beta, \gamma, \theta$ are positive hyperparameters and $\indicator$ is an indicator function. We set $\alpha$ to regulate the establishment of new lightpaths, configure $\beta$ and $\theta$ to prevent extreme imbalances in wavelength occupation on links, and prioritize disjointness between the working route and level-1 routes by setting $\gamma$ and $\chi$. Specifically, when planning a level-1 route for demand $k \in \set{K}'$, we set $\chi = \{e \in \set{E}(r_{k}^{0 0}): \text{route } r_{k}^{e 0} \text{ has not been planned}\}$, i.e., the set of unprotected links in the working route. We set $\chi = \emptyset$ when planning the working route and level-2 routes. Moreover, the number of links traversed is considered in physical arc costs with the cost ``1'' and logical arc costs with the term $|\set{E}(l(a))|$.

\begin{table}
\centering
\caption{Arc costs \label{tab:labeling:cost}}
{\begin{tabular}{ll}
\hline
Arc Type & Cost \\
\hline
Transmitter arc $a = (i', i) \in \set{A}_t$ & $c(a) = \alpha$ \\
Receiver arc $a = (i, i') \in \set{A}_r$ & $c(a) = 0$ \\
Physical arc $a = (i_1, i_2, e) \in \set{A}_p$ & $c(a) = 1 + \beta \left(1 - \frac{|\Lambda(e(a))|}{|\Lambda|}\right)^{\theta} + \gamma \indicator_{\chi} (e(a))$ \\
Logical arc $a = (i_{1}^{\prime}, i_{2}^{\prime}, l) \in \set{A}_l$ & $c(a) = |\set{E}(l(a))| + \gamma |\set{E}(l(a)) \cap \chi|$ \\
\hline
\end{tabular}}
\end{table}

Based on the above cost settings, for any node $j \in \set{N} \cup \set{N}'$, we derive a lower bound on the minimum cost of paths connecting nodes $j$ and $t_{k}^{\prime}$. This estimates the additional cost required to extend a partial path ending at node $j$ to $t_{k}^{\prime}$, and is used to guide the search of routes and accelerate the labeling algorithm. Since arc costs are nonnegative, zero is a trivial lower bound. To obtain a non-trivial, tighter lower bound, we define the minimum hop $h(i_1, i_2)$ as the minimum number of links needed to connect nodes $i_1, i_2 \in \set{N}$. Values of $h(i_1, i_2), i_1, i_2 \in \set{N}$ are computed in the preprocessing stage of Algorithm~\ref{alg:construction:outline}, by setting the length of each link to one and running shortest-path algorithms. The minimum hop extends from the OTN to the LBAG, by defining $h(i_{1}^{\prime}, i_{2}^{\prime}) = h(i_1, i_{2}^{\prime}) = h(i_{1}^{\prime}, i_2) = h(i_1, i_2), i_1, i_2 \in \set{N}$. Proposition~\ref{prop:minhop:LB} shows that $h(j, t_{k}^{\prime})$ is a valid lower bound on the minimum cost of paths connecting nodes $j$ and $t_{k}^{\prime}$.

\begin{proposition}
\label{prop:minhop:LB}
	$h(j, t_{k}^{\prime}), j \in \set{N} \cup \set{N}'$ is a lower bound on the minimum cost of paths connecting nodes $j$ and $t_{k}^{\prime}$.
\end{proposition}

\proof{Proof.}
See Appendix~\ref{proof:minhop:LB}.
\endproof

\subsubsection{Label Definition, Feasibility, and Extension}
\label{sect:labeling:label}


A label $L = (f, j, \bar{c}, \bar{\rho}, d, W, V, U)$ represents a partial or complete path $\rho(L) = (a_1, a_2, \ldots, a_n = f)$ on the LBAG, where arc $a_1$ starts at node $s_{k}^{\prime}$ and  arc $f$ ends at node $j$. The estimated cost $\bar{c} = \sum_{a \in \rho(L)} c(a) + h(j, t_{k}^{\prime})$ is the sum of the current cost and a lower bound (defined in Section~\ref{sect:labeling:cost}) on the additional cost required to reach the terminal node $t_{k}^{\prime}$. $\bar{\rho}$ is the last physical segment, $d$ is the length of this segment, and $W$ is the set of wavelengths unoccupied on the segment. If arc $f$ is not part of the physical layer, i.e., $f \in \set{A}_l \cup \set{A}_t \cup \set{A}_r$, then $\bar{\rho} = \emptyset, d = 0, W = \Lambda$. Otherwise, if $f \in \set{A}_p$, then $\bar{\rho} = (a_{m}, a_{m + 1}, \ldots, a_n)$ where $a_{m}, a_{m + 1}, \ldots, a_n \in \set{A}_p$ and $a_{m - 1} \in \set{A}_t$, $d = \sum_{a \in \bar{\rho}} d(e(a))$, and $W = \bigcap_{a \in \bar{\rho}} \Lambda(e(a))$. $V = \{s_k\} \cup \bigcup_{a \in \rho(L) \cap \left(\set{A}_p \cup \set{A}_l\right)} \set{N}(a)$ is the set of nodes traversed by the physical routing of path $\rho(L)$, where $\set{N}(a)$ is the set of nodes traversed by the physical routing of arc $a$, excluding the first node. $U = \{e(a): a \in \rho(L) \cap \set{A}_p\}$ is the set of links associated with physical arcs traversed by path $\rho(L)$.

\textbf{Feasibility Conditions for Label Extension} When finding a route for demand $k \in \set{K}'$ in scenario $(e_1, e_2) \in \Omega$, the extension of a label $L_1 = (f_{1}, j_1, \bar{c}_1, \bar{\rho}_1, d_1, W_1, V_1, U_1)$ along an arc $f_{2}$ adjacent to node $j_1$ is deemed feasible if and only if all the following practical side constraints and optimality conditions are satisfied:
\begin{itemize}
	\item Simple path constraints in the working scenario: If $(e_1, e_2) = (0, 0)$, then $V_1 \cap \set{N}(f_2) = \emptyset$.
	\item Survivable constraints in failure scenarios: If $(e_1, e_2) \in \Omega^{1} \cup \Omega^{2}$, then $e(f_2) \not\in \{e_1, e_2\}$ when $f_2 \in \set{A}_p$, and $\set{E}(f_2) \cap \{e_1, e_2\} = \emptyset$ when $f_2 \in \set{A}_l$.
	\item Optical reach constraints: If $f_2 \in \set{A}_p$, then $d_1 + d(e(f_2)) \le \bar{d}$.
	\item Wavelength continuity constraints: If $f_2 \in \set{A}_p$, then $W_1 \cap \Lambda(e(f_{2})) \ne \emptyset$.
	\item Bandwidth capacity constraints: If $f_2 \in \set{A}_l$, then $C_{l(f_2)}^{e_1 e_2} + b_k \le B$ when $(e_1, e_2) = (0, 0)$ and $C_{l(f_2)}^{e_1 e_2} + b_k (1 - y_{k l(f_2)}^{0 0}) \le B$ when $(e_1, e_2) \in \Omega^{1} \cup \Omega^{2}$.
	\item Necessary conditions for minimum-cost paths: If $f_2 \in \set{A}_p$, then $e(f_2) \not\in U_1$.
\end{itemize}

\textbf{Label Extension} If label $L_1 = (f_{1}, j_1, \bar{c}_1, \bar{\rho}_1, d_1, W_1, V_1, U_1)$ can be extended along arc $f_{2}$ with head $j_1$ and tail $j_2$, then the concatenation of them yields a label $L_2 = (f_{2}, j_2, \bar{c}_2, \bar{\rho}_2, d_2, W_2, V_2, U_2)$ where (i) $\bar{c}_2 = \bar{c}_1 + c(f_{2}) - h(j_1, t_{k}^{\prime}) + h(j_2, t_{k}^{\prime})$, (ii) if $f_{2} \in \set{A}_l \cup \set{A}_t \cup \set{A}_r$, then $\bar{\rho}_2 = \emptyset, d_2 = 0, W_2 = \Lambda$, (iii) if $f_{2} \in \set{A}_p$, then $\bar{\rho}_2 = \bar{\rho}_1 \cup \{f_{2}\}, d_2 = d_1 + d(e(f_{2})), W_2 = W_1 \cap \Lambda(e(f_{2}))$, (iv) $V_2 = V_1 \cup \set{N}(f_2)$, and (v) if $f_{2} \in \set{A}_p$, then $U_2 = U_1 \cup \{e(f_2)\}$.

The extension procedure follows Algorithm~\ref{alg:labeling}, which iteratively extends labels of the smallest estimated costs to generate new labels and search routes. In line~\ref{alg:labeling:route-generation}, a route $r_{k}^{e_1 e_2}$ is generated by sequentially concatenating lightpaths associated with path $\rho(L_2)$, where new lightpaths might be established (see Appendix~\ref{sect:labeling:lightpath}). New lightpaths are added to $\set{L}'$, and the LBAG $\set{G}$ is updated accordingly. Moreover, we set $y_{k l}^{e_1 e_2} = 1, \forall l \in \set{L}(r_{k}^{e_1 e_2})$.

\begin{algorithm}
\begin{small}
\caption{Labeling Algorithm for Route Planning}
\label{alg:labeling}
\begin{algorithmic}[1]
	\Statex \algorithmicrequire \text{ } Established lightpaths $\set{L}'$, bandwidth consumption $\vect{C}$, partial solution $\vect{y}$, LBAG $\set{G}$, demand $k$, scenario $(e_1, e_2)$
	\Statex \algorithmicensure \text{ } Route $r_{k}^{e_1 e_2}$
	\State Initialize the priority queue of labels as $\set{S} = \{L_0\}$ where $L_0 = (\emptyset, s_{k}^{\prime}, h(s_{k}^{\prime}, t_{k}^{\prime}), \emptyset, 0, \Lambda, \emptyset, \emptyset)$
	\While{$\set{S} \ne \emptyset$}
		\State Find a label $L_1 = (f_1, j_1, \bar{c}_1, \bar{\rho}_1, d_1, W_1, V_1, U_1)$ where $\bar{c}_1 = \min \{\bar{c}(L): L \in \set{S}\}$, i.e., the label in $\set{S}$ with the smallest estimated cost
		\State Let $\set{S} \leftarrow \set{S} \setminus \{L_1\}$
		\ForAll{arc $f_{2} \in \set{G}$ adjacent to node $j_1$}
			\If{$L_1$ can be extended along $f_2$} $\quad$(Including capacity check, Appendix~\ref{sect:bandwidth:check})
				\State Concatenate $L_1$ and $f_{2}$ to form a label $L_2 = (f_2, j_2, \bar{c}_2, \bar{\rho}_2, d_2, W_2, V_2, U_2)$
				\If{$j_2 = t_{k}^{\prime}$}
					\State According to path $\rho(L_2)$, generate route $r_{k}^{e_1 e_2}$ and update $\set{L}'$, $\set{G}$, $\vect{y}$ \label{alg:labeling:route-generation}
					\State Update bandwidth consumption $\vect{C}_{l}, \forall l \in \set{L}(r_{k}^{e_1 e_2})$ $\qquad \quad$ (Appendix~\ref{sect:bandwidth:update})
					\State Terminate the labeling algorithm
				\Else
					\State Apply dominance rules to $L_2$ and labels in $\set{S}$ $\qquad \qquad \quad$ (Section~\ref{sect:labeling:dominance})
					\State Let $\set{S} \leftarrow \set{S} \cup \{L_2\}$ if $L_2$ has not been removed and $W_2 \ne \emptyset$
				\EndIf
			\EndIf
		\EndFor
	\EndWhile
\end{algorithmic}
\end{small}
\end{algorithm}

\subsubsection{Heuristic Dominance Rules}
\label{sect:labeling:dominance}

To reduce the number of labels for exploration, we propose heuristic dominance rules as shown in Table \ref{tab:dominance}. They are relaxations of exact dominance rules so that fewer labels need to be explored at the sacrifice of optimality. Conditions (i) and (ii) require that both labels end at the same node, with the node being on the logical and physical layer, respectively. Condition (iii) stipulates that label $L_1$ must have an equal or smaller cost. Condition (iv) addresses optical reach constraints, applicable to labels ending at a physical node.

For two labels $L_1 = (f_1, j_1, \bar{c}_1, \bar{\rho}_1, d_1, W_1, V_1, U_1)$ and $L_2 = (f_2, j_2, \bar{c}_2, \bar{\rho}_2, d_2, W_2, V_2, U_2)$: If conditions (i) and (iii) specified by Rule 1 are satisfied, we remove $L_2$ from the priority queue $\set{S}$; If conditions (ii)--(iv) specified by Rule 2 are satisfied, we set $W_2 \leftarrow W_2 \setminus W_1$, and further remove $L_2$ from the priority queue $\set{S}$ in case that $W_2$ becomes empty.

Moreover, when condition (v) is satisfied, Rule 2 can be lifted to Rule 3. In this case, $\bigcap_{e \in \set{E}'} \Lambda(e) \ne \emptyset$, i.e., there is a wavelength unoccupied on all links in $\set{E}'$, where $\set{E}' = \{e \in \set{E} \setminus \{e_1, e_2\}: \Lambda(e) \ne \emptyset\}$ is the set of non-failed links with unoccupied wavelengths. Since any path generated by the labeling algorithm must not traverse physical arcs associated with links in $\set{E} \setminus \set{E}'$ and can traverse physical arcs associated with each link $e \in \set{E}'$ at most once, wavelength continuity constraints are always satisfied by assigning the newly established lightpath with a wavelength $\lambda \in \bigcap_{e \in \set{E}'} \Lambda(e)$. Therefore, when conditions (ii)--(v) specified by Rule 3 are satisfied, we directly remove $L_2$ from the priority queue $\set{S}$.

\begin{table}
\centering
\caption{Heuristic dominance rules between two labels $L_1 = (f_1, j_1, \bar{c}_1, \bar{\rho}_1, d_1, W_1, V_1, U_1)$ and $L_2 = (f_2, j_2, \bar{c}_2, \bar{\rho}_2, d_2, W_2, V_2, U_2)$ \label{tab:dominance}}
{\begin{tabular}{lccc}
\hline
 & Rule 1 & Rule 2 & Rule 3 \\
\hline
Condition (i): $j_1 = j_2, j_1 \in \set{N}'$ & \checkmark & & \\
Condition (ii): $j_1 = j_2, j_1 \in \set{N}$ & & \checkmark & \checkmark \\
Condition (iii): $\bar{c}_1 \le \bar{c}_2$ & \checkmark & \checkmark & \checkmark \\
Condition (iv): $d_1 \le d_2$ & & \checkmark & \checkmark \\
Condition (v): $\bigcap_{e \in \set{E} \setminus \{e_1, e_2\}: \Lambda(e) \ne \emptyset} \Lambda(e) \ne \emptyset$ & & & \checkmark \\
\hline
\end{tabular}}
\end{table}

\subsection{Reuse-then-Search Strategy for Route Planning}
\label{sect:reuse-route}

To reduce the number of calls to the labeling algorithm, we adopt a reuse-then-search strategy for route planning. As shown in lines \ref{alg:construction:outline:level-1:route} and \ref{alg:construction:outline:level-2:route} of Algorithm~\ref{alg:construction:outline}, the labeling algorithm is called only when no existing route can be reused. Algorithm~\ref{alg:reuse-route} shows how to identify routes that can be reused. For each pair of terminal nodes $\{s, t\}$, denote $\set{R}_{s t}^{\prime}$ as the set of existing routes that have been generated with endpoints $s$ and $t$. A route $r \in \set{R}_{s t}^{\prime}$ can be reused to protect demand $k$ in the failure scenario $(e_1, e_2) \in \Phi_k \cup \Psi_k \cup \Theta_k$ if and only if: $\{s_k, t_k\} = \{s, t\}$, route $r$ does not traverse any failed link in $\{e_1, e_2\} \cap \set{E}$, and bandwidth capacity constraints are satisfied.

\begin{algorithm}
\begin{small}
\caption{Route Reuse}
\label{alg:reuse-route}
\begin{algorithmic}[1]
	\Statex \algorithmicrequire \text{ } Established lightpaths $\set{L}'$, bandwidth consumption $\vect{C}$, partial solution $\vect{y}$, demand $k$, failure scenario $(e_1, e_2) \in \Phi_k \cup \Psi_k \cup \Theta_k$
	\Statex \algorithmicensure \text{ } Route $r_{k}^{e_1 e_2}$
	\ForAll{route $r \in \set{R}_{s_k t_k}^{\prime}$ \algand $\set{E}(r) \cap \{e_1, e_2\} = \emptyset$}
		\If{Route $r$ satisfy bandwidth capacity constraints} $\qquad \qquad$ (Appendix~\ref{sect:bandwidth:check})
			\State Let $r_{k}^{e_1 e_2} \gets r$, i.e., demand $k$ uses route $r$ in scenario $(e_1, e_2)$
			\State Let $y_{k l}^{e_1 e_2} \gets 1, \forall l \in \set{L}(r_{k}^{e_1 e_2})$
			\State Update bandwidth consumption $\vect{C}_{l}, \forall l \in \set{L}(r_{k}^{e_1 e_2})$ $\qquad \qquad \qquad$(Appendix~\ref{sect:bandwidth:update})
			\State Terminate with successful reuse of route $r_{k}^{e_1 e_2}$
		\EndIf
	\EndFor
	\State Terminate without reuse of routes
\end{algorithmic}
\end{small}
\end{algorithm}

\section{Enhancement by Parallel Computing}
\label{sect:parallel}

For the STG2 instances (with thousands of links) addressed in this paper, millions of scenarios are considered. With the continuous growth of network sizes, the number of scenarios in practical STG2 instances may further increase. To manage the large problem sizes and reduce computing time, we employ multithreaded parallel computing, a common approach for solving large-scale optimization problems \citep{Schryen2020}.
Specifically, as outlined in Algorithm~\ref{alg:parallel:outline}, we adapt the hierarchical constructive heuristic by planning level-2 scenarios in parallel. This approach is justified for two reasons. First, the number of level-2 scenarios ($|\set{E}|^{2} - |\set{E}|$) significantly exceeds that of other scenarios ($1 + |\set{E}|$), making them the primary bottleneck. Second, once the working and level-1 scenarios are planned, the interdependency between level-2 scenarios decreases, allowing for parallel planning.

Nevertheless, two challenges persist. First, different threads may simultaneously establish substitutable lightpaths, potentially increasing the objective value. Second, race conditions (multiple threads access and modify shared data concurrently) necessitate thread-safe data structures and operations, potentially compromising algorithm efficiency. To address these challenges, we adopt a ``master-slave'' framework, executing multiple slave threads and a single master thread alternately. In each iteration, the slave threads plan different scenarios in parallel without establishing new lightpaths. The master thread, capable of establishing new lightpaths, handles scenarios that the slave threads cannot plan (referred to as ``blocked" scenarios), due to insufficient lightpaths. This approach not only resolves the first challenge, but also mitigates race conditions, as $\set{L}'$ and $\set{G}$ become constant when running slave threads in parallel. To further avoid race conditions when accessing and modifying scenario statuses (indicating whether a scenario has been planned successfully, is blocked, or has not been planned), we assign slave threads disjoint pools of scenarios before running them in parallel. We also design stopping criteria to minimize thread idleness.

Algorithm~\ref{alg:parallel:outline} operates as follows. After planning the working scenario and level-1 scenarios (lines~\ref{alg:construction:outline:agg:demand}--\ref{alg:construction:outline:sort:scn:2} of Algorithm~\ref{alg:construction:outline}), it iterates until all level-2 scenarios are planned. Each iteration consists of three main phases:
\begin{itemize}
	\item Phase 1 (lines~\ref{alg:parallel:outline:init:blocked}--\ref{alg:parallel:outline:assign-scenario}). In line~\ref{alg:parallel:outline:init:blocked}, we clear the set of blocked scenarios $\tilde{\Gamma}$, and initialize the bandwidth consumption in blocked scenarios $\tilde{\vect{C}}$. In line~\ref{alg:parallel:outline:assign-scenario}, we initialize the pool of scenarios $\Gamma_{t}$ for slave thread $t = 1, 2, \ldots, M$, as scenarios in $\hat{\Omega}^{2}$ with indices $\{t + i M: i \in \dom{Z}_{+}, i < NS, t + i M \le |\hat{\Omega}^{2}|\}$. Here, $M$ is the number of slave threads, $NS$ is the size limit on a scenario pool, and $\dom{Z}_{+}$ is the set of non-negative integers. This approach balances the workload of slave threads, as scenarios in $\hat{\Omega}^{2}$ are sorted in non-increasing order of the number of remaining demands to be planned.
	\item Phase 2 (line~\ref{alg:parallel:outline:slave}). Plan scenarios in $\Gamma_{1}, \Gamma_{2}, \ldots, \Gamma_{M}$ using slave threads $1, 2, \ldots, M$ in parallel, without establishing new lightpaths. The returned values for each slave thread $t$ include the set of successfully planned scenarios $\bar{\Gamma}_t$, the set of blocked scenarios $\tilde{\Gamma}_t$, and the bandwidth consumption $\tilde{\vect{C}}(t)$ in blocked scenarios. Stopping criteria for slave threads are detailed in Appendix~\ref{sect:parallel:slave}. Particularly, once a thread terminates, all other threads will promptly terminate to minimize thread idleness.
	\item Phase 3 (line~\ref{alg:parallel:outline:master}). Use the master thread to plan blocked scenarios $\tilde{\Gamma}$ that slave threads could not handle, allowing to establish new lightpaths. For details, see Appendix~\ref{sect:parallel:master}.
\end{itemize}

\begin{algorithm}
\begin{small}
\caption{Enhancement by Parallel Computing}
\label{alg:parallel:outline}
\begin{algorithmic}[1]
	\Statex \algorithmicrequire \text{ } $\set{N}, \set{E}, \Lambda, \set{K}, B, \bar{d}$, and hyperparameters $M, NS$
	\Statex \algorithmicensure \text{ } Routes $r_{k}^{e_1 e_2}, k \in \set{K}, (e_1, e_2) \in \Omega$
	\State Plan the working scenario and level-1 scenarios as lines~\ref{alg:construction:outline:agg:demand}--\ref{alg:construction:outline:sort:scn:2} of Algorithm~\ref{alg:construction:outline}
	\State Initialize the set of level-2 scenarios to plan: $\hat{\Omega}^{2} \leftarrow \Omega^{2}$
	\While{$\hat{\Omega}^{2} \ne \emptyset$}
		\State Let blocked scenarios $\tilde{\Gamma} \gets \emptyset$, bandwidth consumption in blocked scenarios $\tilde{\vect{C}} \gets \emptyset$ \label{alg:parallel:outline:init:blocked}
		\State Initialize scenario pools $\Gamma_{1}, \Gamma_{2}, \ldots, \Gamma_{M}$ for threads $1, 2, \ldots, M$ \label{alg:parallel:outline:assign-scenario}
		\State Plan scenarios in $\Gamma_{1}, \Gamma_{2}, \ldots, \Gamma_{M}$ by slave threads $1, 2, \ldots, M$ in parallel: $\bar{\Gamma}_t, \tilde{\Gamma}_t, \tilde{\vect{C}}(t) \gets$ \Call{SlavePlan}{$\Gamma_t$} $\qquad \qquad \quad$ (Slave threads in parallel, see \ref{sect:parallel:slave}) \label{alg:parallel:outline:slave}
		\State Update blocked scenarios and bandwidth consumption in these scenarios: $\tilde{\Gamma} \gets \bigcup_{t = 1}^{M} \tilde{\Gamma}_t, \tilde{\vect{C}} \gets \bigcup_{t = 1}^{M} \tilde{\vect{C}}(t)$
		\State \Call{MasterPlan}{$\tilde{\Gamma}, \tilde{\vect{C}}$} $\qquad \qquad \qquad \qquad \qquad$ (Plan blocked scenarios, see \ref{sect:parallel:master}) \label{alg:parallel:outline:master}
		\State Update unplanned scenarios: $\hat{\Omega}^{2} \gets \hat{\Omega}^{2} \setminus \bigcup_{t = 1}^{M} \bar{\Gamma}_t \setminus \tilde{\Gamma}$ \label{alg:parallel:outline:update}
	\EndWhile
\end{algorithmic}
\end{small}
\end{algorithm}

\section{Computational Experiments}
\label{sect:experiments}

Our computational experiments aim to: (i) demonstrate the advantages of our solution method over the benchmark method of our industry partner (Huawei Technologies), and (ii) examine the effectiveness of algorithm enhancements. The benchmark method is also a constructive heuristic with post-processing to enhance the objective value, and its details are unavailable due to confidentiality. During our study, we developed several other constructive heuristics, including the one that plans demands sequentially with an outer loop iterating over demands and an inner loop over scenarios. Nonetheless, our computational results indicate that these alternative heuristics yield solutions inferior to those from the benchmark method. Furthermore, our hierarchical constructive heuristic significantly outperforms all alternatives in terms of objective value, computing time, and memory usage. Therefore, we report only computational results of the benchmark method and our hierarchical constructive heuristic.

We implement our hierarchical constructive heuristic in C++ and conduct experiments on a server equipped with Intel(R) Xeon(R) Gold 6132 CPU @ 2.60GHz, limiting usage to 16 GB of RAM. The experiments are conducted on 16 large-scale instances provided by our industry partner, featuring 1,000--2,600 nodes, 1,928--3,091 links, 10,020--12,980 demands, and up to over 9 million double-link failure scenarios. In each instance, the number of available wavelengths is $|\Lambda| = 120$ and the bandwidth of each wavelength is $B = 100$ Gbps. For the hyperparameters related to arc costs (Section~\ref{sect:labeling:cost}), we set $\alpha = 32$, $\beta = 32$, $\gamma = 2$, and $\theta = 3$. In the parallel algorithm (Section~\ref{sect:parallel}), we set the size limit of the scenario pool for each slave thread to $NS =$10,000.

Table~\ref{tab:experiments:benchmark:comparison} compares our algorithms with the benchmark method of our industry partner, and demonstrates the effectiveness of the parallel computing enhancement. CH-1 denotes the hierarchical constructive heuristic without the parallel computing enhancement, and CH-M (with $M \ge 2$) denotes the enhanced variant with $M$ slave threads. For each algorithm, column ``Obj'' presents the objective value, and column ``Time'' presents the computing time. Compared to the benchmark method, $Imp_B$ represents the improvement in objective value, and $TR_B$ is the ratio of computing times. Compared to CH-1, $TR_1$ is the ratio of computing times for planning level-2 scenarios (see Appendix~\ref{sect:experiments:parallel} for details).

Table~\ref{tab:experiments:benchmark:comparison} demonstrates that our algorithms outperform the benchmark method in terms of objective value, with the exception of one instance ($|\set{N}| =$2,000). CH-8 achieves an average improvement of 18.5\%. In terms of efficiency, CH-8 uses only 1.5\% of the computing time required by the benchmark method on average, making it 67 times faster. CH-8 averages 14 minutes per instance (with a maximum of 34 minutes), whereas the benchmark method averages 11.3 hours per instance (with a maximum of 21.1 hours). The results also highlight the effectiveness of parallel computing: Compared to CH-1, CH-4 and CH-8 achieve speedups of 3.6 and 6.5, respectively, when planning level-2 scenarios. 

\begin{table}
\begin{small}
\centering
\caption{Computational results of the benchmark method and our algorithms \label{tab:experiments:benchmark:comparison}}
{\begin{tabular}{rrrr|rr|rr|rr|rr}
\hline
\multicolumn{4}{c|}{Instance} & \multicolumn{2}{c|}{Benchmark} & \multicolumn{2}{c|}{CH-1} & \multicolumn{2}{c|}{CH-4} & \multicolumn{2}{c}{CH-8} \\
\hline
$|\set{N}|$ & $|\set{E}|$ & $|\set{K}|$ & $\bar{d}$ (km) & Obj  & Time (s) & Obj  & Time (s) & Obj  & Time (s) & Obj  & Time (s)\\
\hline
1000  & 1928  & 10020  & 360   & 7151  & 51051 & 5959   & 1722.2  & 5862   & 689.7  & 5882   & 393.6  \\
1000  & 1928  & 10020  & 600   & 7222  & 37080 & 5883   & 1719.8  & 5863   & 695.2  & 5870   & 394.1  \\
1000  & 1928  & 10020  & 900   & 6870  & 34603 & 5179   & 1296.2  & 5176   & 559.5  & 5131   & 305.2  \\
1000  & 1928  & 10020  & 1500  & 6728  & 34305 & 4763   & 1100.9  & 4805   & 460.8  & 4792   & 262.4  \\
1200  & 2365  & 10390  & 600   & 7927  & 37137 & 6675   & 2703.9  & 6589   & 1071.5 & 6605   & 564.7  \\
1200  & 2365  & 10390  & 900   & 7615  & 28259 & 5903   & 2031.4  & 5856   & 863.4  & 5874   & 458.0  \\
1200  & 2365  & 10390  & 1500  & 7685  & 28017 & 5513   & 1662.1  & 5486   & 720.4  & 5514   & 394.8  \\
1600  & 3091  & 10070  & 600   & 8260  & 73270 & 6700   & 4767.2  & 6673   & 1713.8 & 6658   & 837.4  \\
1600  & 3091  & 10070  & 900   & 7473  & 75978 & 5804   & 5313.4  & 5750   & 1076.8 & 5779   & 663.8  \\
2000  & 3083  & 12980  & 600   & 7596  & 20845 & 8543   & 14495.6 & 8496   & 2770.0 & 8540   & 1695.0 \\
2200  & 3052  & 11054  & 600   & 10405 & 37137 & 8925   & 17661.9 & 8839   & 4107.7 & 8882   & 2031.5 \\
2200  & 3052  & 11054  & 900   & 8540  & 44990 & 6858   & 10624.1 & 6850   & 2054.7 & 6856   & 1284.1 \\
2200  & 3052  & 11054  & 1500  & 8208  & 38184 & 5591   & 4911.3  & 5545   & 1535.1 & 5541   & 949.1  \\
2600  & 3027  & 10930  & 600   & 10151 & 68144 & 8924   & 9671.7  & 8822   & 1992.4 & 8908   & 1254.4 \\
2600  & 3027  & 10930  & 900   & 7729  & 27941 & 6588   & 5363.5  & 6559   & 1302.2 & 6573   & 817.6  \\
2600  & 3027  & 10930  & 1500  & 5749  & 18637 & 4990   & 3657.8  & 4953   & 978.9  & 4942   & 607.2  \\
\hline
\multicolumn{6}{c|}{Average $Imp_B$ and $TR_B$} & 18.1\% & 0.074 & 18.7\% & 0.026 & 18.5\% & 0.015 \\
\hline
\multicolumn{6}{c|}{Average $TR_1$} & \multicolumn{2}{c|}{1.0} & \multicolumn{2}{c|}{3.6} & \multicolumn{2}{c}{6.5}\\
\hline
\end{tabular}}
\end{small}
\end{table}

Table~\ref{tab:experiments:effectiveness:symbols} presents the notations for the hierarchical constructive heuristic with various enhancements. All these algorithm variants are enhanced with parallel computing (Section~\ref{sect:parallel}) using 8 slave threads, with CH equivalent to CH-8 for brevity. Disabeling dominance rule 1 results in failure to solve the smallest-scale instances (1,000 nodes) within a 128 GB memory limit. In addition to examining the effectiveness of demand aggregation and dominance rules, we also compare the results with and without using the bound on path costs in the labeling algorithm, i.e., calculating the estimated cost of a label associated with path $\rho = (a_1, a_2, \ldots, a_n)$ ending at node $j$ with $\bar{c} = \sum_{a \in \rho} c(a) + h(j, t_{k}^{\prime})$ and $\bar{c} = \sum_{a \in \rho} c(a)$, respectively. These algorithm variants are compared to CH.

\begin{table}
\centering
\caption{Hierarchical Constructive Heuristic with Various Enhancements \label{tab:experiments:effectiveness:symbols}}
{\begin{tabular}{lcccc}
\hline
Notation & Demand Aggregation & \makecell{Use Bound \\ on Path Costs} & Route Reuse & Dominance Rules 2 and 3 \\
\hline
CH      & \checkmark  & \checkmark & \checkmark & \checkmark \\
CH-D    &   & \checkmark & \checkmark & \checkmark \\
CH-B    & \checkmark &   & \checkmark & \checkmark \\
CH-R    & \checkmark & \checkmark &   & \checkmark \\
CH-R-DP & \checkmark & \checkmark &   &   \\
\hline
\end{tabular}}
\end{table}

By comparing the computational results of CH-D and CH, which are shown in Table~ \ref{tab:experiments:effectiveness:aggregation} of Appendix~\ref{app:sec:enhance}, we evaluate the effectiveness of demand aggregation. The results show that without demand aggregation, the objective value worsens 7.6\% on average and computing time lengthens 4.3 times on average. Moreover, by comparing the computational results of CH-R-DP, CH-R, CH-B, and CH, which are shown in Table~\ref{tab:experiments:effectiveness:routing} of Appendix~\ref{app:sec:enhance}, we evaluate the performance of different routing techniques. 
The results show that the computing time decreases significantly (by 39.4 times on average) with route reuse adopted, and it decreases further with dominance rules 2 and 3 adopted. Additionally, the results also show that CH-R-DP can solve only three instances within the memory limit, highlighting the importance of dominance rules 2 and 3 in reducing memory requirements. Furthermore, incorporating lower bounds on path costs provides a substantial improvement in computing time, with CH using only about half the computing time of CH-B.

\section{Conclusion}
\label{sect:conclusion}

This paper addresses large-scale STG2 with practical constraints, essential for minimizing network costs and ensuring reliable telecommunications  services. To effectively manage the interdependency between scenarios imposed by various constraints, we have developed a hierarchical constructive heuristic. It involves planning the working scenario first, followed by all level-1 scenarios, and finally all level-2 scenarios. This hierarchical sequence offers several advantages: It efficiently handles consistent routing constraints through bandwidth reservation, allows bandwidth capacity constraints to be verified by querying bandwidth consumption only in the scenario being currently planned, facilitates bandwidth updates requiring fewer operations and less memory, and mitigates race conditions during parallel planning of level-2 scenarios.

To address the challenges of large problem sizes, we propose several strategies. To reduce the demand set, we aggregate demands using bin-packing solutions. For efficient route searching in large-scale networks, the labeling algorithm is enhanced by incorporating lower bounds on path costs and heuristic dominance rules, alongside a reuse-then-search routing strategy to reduce calls of the labeling algorithm. To manage the large number of level-2 scenarios, a multithreaded parallel algorithm is implemented, with mechanisms designed to optimize the objective value, mitigate race conditions, and reduce thread idleness.

The computational results indicate that our method efficiently solves large-scale STG2 instances, comprising over 2,000 nodes, 3,000 links, 10,000 demands, and 9 million failure scenarios, within an hour of computing time and 16 GB of memory. This performance significantly surpasses the scales addressed in existing traffic grooming studies. Compared to the benchmark method from our industry partner, our approach yields an average improvement of 18.5\% in objective value and is 67 times faster in terms of computing time. Computational results also demonstrate that algorithm enhancements proposed in this paper are effective. These findings underscore the superiority of our method in objective value optimization, scalability, and computational efficiency.

In future work, we intend to explore additional techniques, such as distributed computing and post-processing, to enhance our hierarchical constructive heuristic. We also aim to develop new efficient constructive methods for addressing large-scale STG2. Additionally, we plan to study relaxations for calculating lower bounds on the optimal objective value. These relaxations will help assess the quality of the incumbent solution and facilitate the derivation of near-optimal solutions. Furthermore, there is considerable interest in adapting the strategies proposed in this paper to optimize other types of large-scale networks, such as power networks and logistics networks, while taking failures into account.


%
%
%



\bibliography{references}




\newpage
\renewcommand{\thepage}{A.\arabic{page}}
\setcounter{page}{1}
\renewcommand{\theequation}{A.\arabic{equation}}

\renewcommand{\thetable}{A.\arabic{table}}
\renewcommand{\thealgorithm}{A.\arabic{algorithm}}

\setcounter{equation}{0}
{\large \bf \noindent Appendices}

\setcounter{table}{0}
\setcounter{algorithm}{0}

\begin{appendices}
\section{Proofs of Propositions}

\subsection{Proof of Proposition~\ref{prop:minhop:LB}}
\label{proof:minhop:LB}
\proof{Proof.}
Obviously, $h(i_1, i_2), i_1, i_2 \in \set{N}$ satisfy triangle inequalities, i.e., $h(i_1, i_2) + h(i_2, i_3) \ge h(i_1, i_3), i_1, i_2, i_3 \in \set{N}$. We define surrogate arc costs and show that they are nonnegative, i.e., $\bar{c}(a) \ge 0, a \in \set{A}_p \cup \set{A}_l \cup \set{A}_t \cup \set{A}_r$:
\begin{itemize}
	\item For transmitter arc $a = (i', i) \in \set{A}_t$, the surrogate cost $\bar{c}(a) = c(a) - h(i', t_{k}^{\prime}) + h(i, t_{k}^{\prime}) = c(a) = \alpha > 0$, since $h(i', t_{k}^{\prime}) = h(i, t_{k}^{\prime})$.
	\item For receiver arc $a = (i, i') \in \set{A}_r$, the surrogate cost $\bar{c}(a) = c(a) - h(i, t_{k}^{\prime}) + h(i', t_{k}^{\prime}) = c(a) = 0$, since $h(i, t_{k}^{\prime}) = h(i', t_{k}^{\prime})$.
	\item For physical arc $a = (i_1, i_2, e) \in \set{A}_p$, the surrogate cost $\bar{c}(a) = c(a) - h(i_1, t_{k}^{\prime}) + h(i_2, t_{k}^{\prime}) = c(a) - h(i_1, t_{k}) + h(i_2, t_{k}) \ge c(a) - h(i_1, i_2) = c(a) - 1 \ge 0$, where $h(i_1, i_2) + h(i_2, t_{k}) \ge h(i_1, t_{k})$ due to triangle inequalities and $h(i_1, i_2) = 1$ since a single link $e$ connects nodes $i_1$ and $i_2$.
	\item For logical arc $a = (i_{1}^{\prime}, i_{2}^{\prime}, l) \in \set{A}_l$, the surrogate cost $\bar{c}(a) = c(a) - h(i_{1}^{\prime}, t_{k}^{\prime}) + h(i_{2}^{\prime}, t_{k}^{\prime}) \ge |\set{E}(l(a))| - h(i_1, t_{k}) + h(i_2, t_{k}) \ge |\set{E}(l(a))| - h(i_1, i_2) \ge 0$, where $h(i_1, i_2) + h(i_2, t_{k}) \ge h(i_1, t_{k})$ due to triangle inequalities and $|\set{E}(l(a))| \ge h(i_1, i_2)$ since lightpath $l(a)$ corresponds to a path connecting nodes $i$ and $j$.
\end{itemize}

For the shortest path $(a_1, a_2, \ldots, a_n)$ connecting nodes $j \in \set{N} \cup \set{N}'$ and $t_{k}^{\prime}, k \in \set{K}'$ on the LBAG, assume that the endpoints of arcs $a_1, a_2, \ldots, a_n$ are $j_1 = j, j_2, \ldots, j_{n}, j_{n + 1} = t_{k}^{\prime}$. Then we have
\begin{subequations}
\begin{eqnarray*}
	&& \sum_{m = 1}^{n} \bar{c}(a_m) = \sum_{m = 1}^{n} \left[c(a_m) - h(j_m, t_{k}^{\prime}) + h(j_{m + 1}, t_{k}^{\prime})\right] \\
	= && \sum_{m = 1}^{n} c(a_m) - h(j_1, t_{k}^{\prime}) + h(j_{n + 1}, t_{k}^{\prime}) = \sum_{m = 1}^{n} c(a_m) - h(j, t_{k}^{\prime}) + h(t_{k}^{\prime}, t_{k}^{\prime}) \\
	= && \sum_{m = 1}^{n} c(a_m) - h(j, t_{k}^{\prime}) \ge 0,
\end{eqnarray*}
\end{subequations}
where $\sum_{m = 1}^{n} c(a_m)$ is the minimum cost of paths connecting nodes $j$ and $t_{k}^{\prime}$, and the inequality follows from that $\bar{c}(a) \ge 0, a \in \set{A}_p \cup \set{A}_l \cup \set{A}_t \cup \set{A}_r$. This implies that $h(j, t_{k}^{\prime})$ is a lower bound on the minimum cost of paths connecting nodes $j$ and $t_{k}^{\prime}$.
\endproof

\subsection{Proof of Proposition~\ref{prop:relation:bandwidth}}
\label{proof:relation:bandwidth}
\proof{Proof.}

Part (i). Algorithm~\ref{alg:construction:outline} starts route planning at line~\ref{alg:construction:outline:working:begin}, thus $C_{l}^{e_1 e_2} = 0, \forall l \in \set{L}, (e_1, e_2) \in \Omega$ before line~\ref{alg:construction:outline:working:begin}, i.e., relations~\eqref{eqn:bandwidth:working:equal:EC} hold before planning the working routes for demands.
\begin{equation}
\label{eqn:bandwidth:working:equal:EC}
	C_{l}^{e_1 e_2} = C_{l}^{0 0}, \forall l \in \set{L}, (e_1, e_2) \in \Omega^{1} \cup \Omega^{2}.
\end{equation}

Throughout lines~\ref{alg:construction:outline:working:begin}--\ref{alg:construction:outline:working:end}, Algorithm~\ref{alg:construction:outline} iteratively plans the working route for each demand. When assigning route $r_{k}^{0 0}$ to demand $k \in \set{K}'$ in the working scenario, the bandwidth consumption on lightpaths in the working scenario should be updated according to operations~\eqref{eqn:bandwidth:update:working:working:EC}, and the bandwidth consumption on lightpaths in failure scenarios should be updated according to operations~\eqref{eqn:bandwidth:update:working:failure:EC} due to bandwidth reservation constraints.
\begin{subequations}
\begin{eqnarray}
	&& C_{l}^{0 0} \gets C_{l}^{0 0} + b_k, \qquad \forall l \in \set{L}(r_{k}^{0 0}), \label{eqn:bandwidth:update:working:working:EC} \\
	&& C_{l}^{e_1 e_2} \gets C_{l}^{e_1 e_2} + b_k, \qquad \forall l \in \set{L}(r_{k}^{0 0}), (e_1, e_2) \in \Omega^{1} \cup \Omega^{2}. \label{eqn:bandwidth:update:working:failure:EC}
\end{eqnarray}
\end{subequations}
Operations~\eqref{eqn:bandwidth:update:working:working:EC} and \eqref{eqn:bandwidth:update:working:failure:EC} indicate that planning working routes for demands (lines~\ref{alg:construction:outline:working:begin}--\ref{alg:construction:outline:working:end}) does not change the fulfillment of relations~\eqref{eqn:bandwidth:working:equal:EC}, which, together with the fact that relations~\eqref{eqn:bandwidth:working:equal:EC} hold before line~\ref{alg:construction:outline:working:begin}, implies that relations~\eqref{eqn:bandwidth:working:equal:EC} hold throughout lines~\ref{alg:construction:outline:working:begin}--\ref{alg:construction:outline:working:end}.

Part (ii). Operations in line~\ref{alg:construction:outline:sort:scn:1}
do not change the fulfillment of relations~\eqref{eqn:bandwidth:working:equal:EC}, which, together with the fact that relations~\eqref{eqn:bandwidth:working:equal:EC} hold throughout lines~\ref{alg:construction:outline:working:begin}--\ref{alg:construction:outline:working:end}, implies that relations~\eqref{eqn:bandwidth:level-2:le:level-1:EC} hold before line~\ref{alg:construction:outline:level-1:begin}.
\begin{equation}
\label{eqn:bandwidth:level-2:le:level-1:EC}
	C_{l}^{e_1 e_2} \le C_{l}^{e_1 0}, \forall l \in \set{L}, (e_1, e_2) \in \Omega^{2}.
\end{equation}

Throughout lines~\ref{alg:construction:outline:level-1:begin}--\ref{alg:construction:outline:level-1:end}, Algorithm~\ref{alg:construction:outline} iteratively plans level-1 routes for demands. When assigning route $r_{k}^{e_1 0}$ to demand $k \in \set{K}'$ in a level-1 scenario $(e_1, 0) \in \Omega^{1}$, the bandwidth consumption on lightpaths in the level-1 scenario $(e_1, 0)$ should be updated according to operations~\eqref{eqn:bandwidth:update:level1:level1:EC}, and the bandwidth consumption on lightpaths in level-2 scenarios $\bar{\Theta}_{k}(e_1) = \{(e_1, e_2) \in \Omega^{2}: e_2 \in \set{E} \setminus \set{E}(r_{k}^{e_1 0})\}$ should be updated according to operations~\eqref{eqn:bandwidth:update:level1:level2:EC} due to consistent routing constraints which require that demand $k$ should still use route $r_{k}^{e_1 0}$ in scenarios $\bar{\Theta}_{k}(e_1)$.
\begin{subequations}
\begin{eqnarray}
	&& C_{l}^{e_1 0} \gets C_{l}^{e_1 0} + b_k (1 - y_{k l}^{0 0}), \qquad \forall l \in \set{L}(r_{k}^{e_1 0}), \label{eqn:bandwidth:update:level1:level1:EC} \\
	&& C_{l}^{e_1 e_2} \gets C_{l}^{e_1 e_2} + b_k (1 - y_{k l}^{0 0}), \qquad \forall l \in \set{L}(r_{k}^{e_1 0}), (e_1, e_2) \in \bar{\Theta}_{k}(e_1). \label{eqn:bandwidth:update:level1:level2:EC}
\end{eqnarray}
\end{subequations}
Operations~\eqref{eqn:bandwidth:update:level1:level1:EC} and \eqref{eqn:bandwidth:update:level1:level2:EC} indicate that planning level-1 routes for demands does not change the fulfillment of relations~\eqref{eqn:bandwidth:level-2:le:level-1:EC}, which, together with the fact that relations~\eqref{eqn:bandwidth:level-2:le:level-1:EC} hold before line~\ref{alg:construction:outline:level-1:begin}, implies that \eqref{eqn:bandwidth:working:equal:EC} hold throughout lines~\ref{alg:construction:outline:level-1:begin}--\ref{alg:construction:outline:level-1:end}.
\endproof

\section{Mathematical Formulation for STG2}
\label{sect:formulation}

For each demand $k \in \set{K}$, we define $\set{R}_k \subseteq \set{R}$ as the subset of routes with endpoints the same with terminal nodes $s_k$ and $t_k$, and $\bar{\set{R}}_k \subseteq \set{R}_k$ as the set of elementary routes (traverse each node at most once) in $\set{R}_k$.

Let $x_l, l \in \set{L}$ be a binary variable equal to one if and only if lightpath $l$ is established. Let $y_{k l}^{e_1 e_2}, k \in \set{K}, l \in \set{L}, (e_1, e_2) \in \Omega: e_1, e_2 \not\in \set{E}(l)$ be a binary variable equal to one if and only if lightpath $l$ belongs to the route used by demand $k$ in scenario $(e_1, e_2)$. Let $z_{k r}^{0 0}, k \in \set{K}, r \in \bar{\set{R}}_k$ be a binary variable equal to one if and only if route $r$ is the working route of demand $k$. Let $z_{k r}^{e_1 e_2}, k \in \set{K}, r \in \set{R}_k, (e_1, e_2) \in \Omega^{1} \cup \Omega^{2}: e_1, e_2 \not\in \set{E}(r)$ be a binary variable equal to one if and only if route $r$ is the backup route of demand $k$ in scenario $(e_1, e_2)$. Then the mathematical formulation can be represented as follows.

\begin{subequations}
\begin{eqnarray}
	&& \min \sum_{l \in \set{L}} x_l, \label{eqn:TG2:Obj} \\
	{\rm s.t.} \quad && \sum_{l \in \set{L}: \lambda(l) = \lambda, e \in \set{E}(l)} x_l \le 1, \quad \forall \lambda \in \Lambda, e \in \set{E}, \label{eqn:TG2:WavelengthAssignment} \\
	&& \sum_{r \in \bar{\set{R}}_{k}} z_{k r}^{0 0} = 1, \quad \forall k \in \set{K}, \label{eqn:TG2:WorkRoute} \\
	&& \sum_{r \in \set{R}_{k}: e_1, e_2 \not\in \set{E}(r)} z_{k r}^{e_1 e_2} = 1, \quad \forall k \in \set{K}, (e_1, e_2) \in \Omega^{1} \cup \Omega^{2}, \label{eqn:TG2:BackupRoute} \\
	&& z_{k r}^{0 0} - z_{k r}^{e_1 0} \le 0, \quad \forall k \in \set{K}, r \in \bar{\set{R}}_{k}, (e_1, 0) \in \Omega^{1}: e_1 \not\in \set{E}(r), \label{eqn:TG2:StillUseWorkRoute} \\
	&& z_{k r}^{e_1 0} - z_{k r}^{e_1 e_2} \le 0, \qquad \forall k \in \set{K}, r \in \set{R}_{k}, (e_1, e_2) \in \Omega^{2}: e_1, e_2 \not\in \set{E}(r), \label{eqn:TG2:StillUseBackup1Route} \\
	&& y_{k l}^{0 0} - \sum_{r \in \bar{\set{R}}_{k}: l \in \set{L}(r)} z_{k r}^{0 0} = 0, \quad \forall k \in \set{K}, l \in \set{L}, \label{eqn:TG2:LinkYAndZWork} \\
	&& y_{k l}^{e_1 e_2} - \sum_{r \in \set{R}_{k}: l \in \set{L}(r), e_1, e_2 \not\in \set{E}(r)} z_{k r}^{e_1 e_2} = 0, \nonumber \\
	&& \qquad \qquad \qquad \qquad \qquad \forall k \in \set{K}, l \in \set{L}, (e_1, e_2) \in \Omega^{1} \cup \Omega^{2}: e_1, e_2 \not\in \set{E}(l), \label{eqn:TG2:LinkYAndZFailure} \\
	&& \sum_{k \in \set{K}} b_k y_{k l}^{0 0} - B x_l \le 0, \quad \forall l \in \set{L}, \label{eqn:TG2:CapacityWork} \\
	&& \sum_{k \in \set{K}} b_k \max\{y_{k l}^{0 0}, y_{k l}^{e_1 e_2}\} - B x_l \le 0, \quad \forall l \in \set{L}, (e_1, e_2) \in \Omega^{1} \cup \Omega^{2}: e_1, e_2 \not\in \set{E}(l), \label{eqn:TG2:CapacityFailure} \\
	&& x_l \in \{0, 1\}, \quad \forall l \in \set{L}, \label{eqn:TG2:Domain:X} \\
	&& y_{k l}^{e_1 e_2} \in \{0, 1\}, \quad \forall k \in \set{K}, l \in \set{L}, (e_1, e_2) \in \Omega: e_1, e_2 \not\in \set{E}(l), \label{eqn:TG2:Domain:Y} \\
	&& z_{k r}^{0 0} \in \{0, 1\}, \quad \forall k \in \set{K}, r \in \bar{\set{R}}_k, \label{eqn:TG2:Domain:ZWork} \\
	&& z_{k r}^{e_1 e_2} \in \{0, 1\}, \quad \forall k \in \set{K}, r \in \set{R}_k, (e_1, e_2) \in \Omega^{1} \cup \Omega^{2}: e_1, e_2 \not\in \set{E}(r). \label{eqn:TG2:Domain:ZFailure}
\end{eqnarray}
\end{subequations}

The objective function \eqref{eqn:TG2:Obj} is to minimize the number of established lightpaths. Wavelength assignment constraints \eqref{eqn:TG2:WavelengthAssignment} require that a wavelength on a link can be assigned to at most one established lightpath. Constraints \eqref{eqn:TG2:WorkRoute} ensure that each demand is assigned an elementary working route, while constraints \eqref{eqn:TG2:BackupRoute} require that each demand is assigned a backup route in each failure scenario, avoiding failed links. Consistent routing constraints \eqref{eqn:TG2:StillUseWorkRoute} impose that the working route should still be used as the backup route in a level-1 scenario $(e_1, 0) \in \Omega^{1}$ if it does not traverse the failed link $e_1$. Consistent routing constraints \eqref{eqn:TG2:StillUseBackup1Route} require that the backup route in a level-1 scenario $(e_1, 0) \in \Omega^{1}$ should still be used as the backup route in a level-2 scenario $(e_1, e_2) \in \Omega^{2}$ if it does not traverse the second failed link $e_2$. Constraints \eqref{eqn:TG2:LinkYAndZWork} and \eqref{eqn:TG2:LinkYAndZFailure} link variables $\vect{y}$ and $\vect{z}$. Bandwidth capacity constraints are described with \eqref{eqn:TG2:CapacityWork} and \eqref{eqn:TG2:CapacityFailure}. Constraints \eqref{eqn:TG2:Domain:X}-\eqref{eqn:TG2:Domain:ZFailure} specify domains of variables.

Note that wavelength continuity constraints and optical reach constraints are respected according to the definition of $\set{L}$. Furthermore, the above formulation can be converted into an equivalent ILP by replacing constraints \eqref{eqn:TG2:CapacityFailure} with constraints \eqref{eqn:TG2:LinkWAndYWork:ILP}--\eqref{eqn:TG2:Domain:WFailure:ILP}.
\begin{subequations}
\begin{eqnarray}
	&& w_{k l}^{e_1 e_2} \ge y_{k l}^{0 0}, \quad \forall l \in \set{L}, (e_1, e_2) \in \Omega^{1} \cup \Omega^{2}: e_1, e_2 \not\in \set{E}(l), \label{eqn:TG2:LinkWAndYWork:ILP} \\
	&& w_{k l}^{e_1 e_2} \ge y_{k l}^{e_1 e_2}, \quad \forall l \in \set{L}, (e_1, e_2) \in \Omega^{1} \cup \Omega^{2}: e_1, e_2 \not\in \set{E}(l), \label{eqn:TG2:LinkWAndYFailure:ILP} \\
	&& \sum_{k \in \set{K}} b_k w_{k l}^{e_1 e_2} - B x_l \le 0, \quad \forall l \in \set{L}, (e_1, e_2) \in \Omega^{1} \cup \Omega^{2}: e_1, e_2 \not\in \set{E}(l), \label{eqn:TG2:CapacityFailure:ILP} \\
	&& w_{k l}^{e_1 e_2} \in \{0, 1\}, \quad \forall l \in \set{L}, (e_1, e_2) \in \Omega^{1} \cup \Omega^{2}: e_1, e_2 \not\in \set{E}(l), \label{eqn:TG2:Domain:WFailure:ILP}
\end{eqnarray}
\end{subequations}
where $w_{k l}^{e_1 e_2}$ is a binary variable equal to one if the bandwidth of $b_k$ on lightpath $l$ is consumed ($y_{k l}^{e_1 e_2} = 1$) or reserved ($y_{k l}^{0 0} = 1$ and $y_{k l}^{e_1 e_2} = 0$) for demand $k$ in failure scenario $(e_1, e_2) \in \Omega^{1} \cup \Omega^{2}$.

\section{More Details on Hierarchical Constructive Heuristic}

\subsection{Establishing New Lightpaths}
\label{sect:labeling:lightpath}

In line~\ref{alg:labeling:route-generation} of Algorithm~\ref{alg:labeling}, the route $r_{k}^{e_1 e_2}$ is generated by sequentially concatenating lightpaths associated with path $\rho(L_2)$, where new lightpaths might be established. Specifically, if $\rho(L_2) = (\ldots, a_m, a_{m + 1}, \ldots, a_n, a_{n + 1}, \ldots)$, where $a_{m - 1} \in \set{A}_t$, $a_m, a_{m + 1}, \ldots, a_n \in \set{A}_p$, and $a_{n + 1} \in \set{A}_r$, then a new lightpath $(\lambda, (e(a_m), e(a_{m + 1}), \ldots, e(a_n)))$ is established, where $\lambda = \min \bigcap_{i = m, m + 1, \ldots, n} \Lambda(e(a_i))$. This means choosing the wavelength with the smallest index among those unoccupied on links $e(a_m), e(a_{m + 1}), \ldots, e(a_n)$, to increase the likelihood that there is a wavelength unoccupied on most links, thereby increasing the potential to apply dominance rule 3 (see Section~\ref{sect:labeling:dominance}).

\subsection{Efficient Check of Bandwidth Capacity Constraints}
\label{sect:bandwidth:check}

When finding a route for a demand in a scenario, both the labeling algorithm and the route reuse algorithm involve checking bandwidth capacity constraints. For this, we first introduce a direct check approach, which is time-consuming, and then present a more efficient partial check approach, which significantly reduces the number of inequalities to check.

In the direct check approach shown in Table~\ref{tab:bandwidth:check}, when assigning a route to demand $k \in \set{K}'$ in the working scenario (level-1 scenario, level-2 scenario, respectively), inequalities~\eqref{eqn:bandwidth:check:working:working} and \eqref{eqn:bandwidth:check:working:failure} (\eqref{eqn:bandwidth:check:level1:level1} and \eqref{eqn:bandwidth:check:level1:level2}, \eqref{eqn:bandwidth:check:level2}, respectively) need to be checked to determine whether bandwidth capacity constraints are satisfied, where each $C_{l}^{e_1 e_2}$ is defined in Section~\ref{sect:labeling}. Inequalities~\eqref{eqn:bandwidth:check:working:failure} are considered for bandwidth reservation constraints, and inequalities~\eqref{eqn:bandwidth:check:level1:level2} are considered due to consistent routing constraints (i.e., demand $k$ should continue using route $r_{k}^{e_1 0}$ in scenarios $\bar{\Theta}_{k}(e_1)$).
\begin{subequations}
	\begin{eqnarray}
		&& C_{l}^{0 0} + b_k \le B, \qquad \forall l \in \set{L}(r_{k}^{0 0}). \label{eqn:bandwidth:check:working:working} \\
		&& C_{l}^{e_1 e_2} + b_k \le B, \qquad \forall l \in \set{L}(r_{k}^{0 0}), (e_1, e_2) \in \Omega^{1} \cup \Omega^{2}. \label{eqn:bandwidth:check:working:failure} \\
		&& C_{l}^{e_1 0} + b_k (1 - y_{k l}^{0 0}) \le B, \qquad \forall l \in \set{L}(r_{k}^{e_1 0}). \label{eqn:bandwidth:check:level1:level1} \\
		&& C_{l}^{e_1 e_2} + b_k (1 - y_{k l}^{0 0}) \le B, \qquad \forall l \in \set{L}(r_{k}^{e_1 0}), (e_1, e_2) \in \bar{\Theta}_{k}(e_1). \label{eqn:bandwidth:check:level1:level2} \\
		&& C_{l}^{e_1 e_2} + b_k (1 - y_{k l}^{0 0}) \le B, \qquad \forall l \in \set{L}(r_{k}^{e_1 e_2}). \label{eqn:bandwidth:check:level2}
	\end{eqnarray}
\end{subequations}
The large cardinalities of $\set{E}$ and $\set{E} \setminus \set{E}(r_{k}^{e_1 0})$ in large-scale STG2 instances lead to a large number of inequalities~\eqref{eqn:bandwidth:check:working:failure} and \eqref{eqn:bandwidth:check:level1:level2} ($|\set{L}(r_{k}^{0 0})| |\set{E}|^{2}$ and $|\set{L}(r_{k}^{e_1 0})| |\set{E} \setminus \set{E}(r_{k}^{e_1 0})|$, respectively).

Proposition~\ref{prop:relation:bandwidth} established below enables us to use the partial check approach specified in Table~\ref{tab:bandwidth:check}, eliminating the need to verify inequalities~\eqref{eqn:bandwidth:check:working:failure} and \eqref{eqn:bandwidth:check:level1:level2}.

\begin{proposition}
	\label{prop:relation:bandwidth}
	When executing lines~\ref{alg:construction:outline:working:begin}--\ref{alg:construction:outline:working:end} of Algorithm~\ref{alg:construction:outline}, we have 
	\begin{equation}
		\label{eqn:relation:bandwidth:working}
		C_{l}^{e_1 e_2} = C_{l}^{0 0}, \forall l \in \set{L}, (e_1, e_2) \in \Omega^{1} \cup \Omega^{2}.
	\end{equation}
	
	When executing lines~\ref{alg:construction:outline:level-1:begin}--\ref{alg:construction:outline:level-1:end} of Algorithm~\ref{alg:construction:outline}, we have 
	\begin{equation}
		\label{eqn:relation:bandwidth:level-1}
		C_{l}^{e_1 e_2} \le C_{l}^{e_1 0}, \forall l \in \set{L}, (e_1, e_2) \in \Omega^{2}.
	\end{equation}
\end{proposition}

\proof{Proof.}
See Appendix~\ref{proof:relation:bandwidth}.
\endproof

\begin{table}
\centering
\caption{Inequalities for checking bandwidth capacity constraints \label{tab:bandwidth:check}}
	{\begin{tabular}{lll}
			\hline
			& Direct Check & Partial Check \\
			\hline
			Assign $r_{k}^{0 0}$ in scenario $(0, 0)$ & \eqref{eqn:bandwidth:check:working:working} and \eqref{eqn:bandwidth:check:working:failure} & \eqref{eqn:bandwidth:check:working:working} \\
			Assign $r_{k}^{e_1 0}$ in scenario $(e_1, 0) \in \Phi_k$ & \eqref{eqn:bandwidth:check:level1:level1} and \eqref{eqn:bandwidth:check:level1:level2} & \eqref{eqn:bandwidth:check:level1:level1} \\
			Assign $r_{k}^{e_1 e_2}$ in scenario $(e_1, e_2) \in \Psi_k \cup \Theta_k$ & \eqref{eqn:bandwidth:check:level2} & \eqref{eqn:bandwidth:check:level2} \\
			\hline
	\end{tabular}}
\end{table}

The correctness of the partial check approach is based on the following reasons. When planning working routes (lines~\ref{alg:construction:outline:working:begin}--\ref{alg:construction:outline:working:end} of Algorithm~\ref{alg:construction:outline}), satisfying inequalities~\eqref{eqn:bandwidth:check:working:working} ensures the fulfillment of inequalities~\eqref{eqn:bandwidth:check:working:failure}, as per relations~\eqref{eqn:relation:bandwidth:working}. Similarly, when planning level-1 routes (lines~\ref{alg:construction:outline:level-1:begin}--\ref{alg:construction:outline:level-1:end} of Algorithm~\ref{alg:construction:outline}), satisfying inequalities~\eqref{eqn:bandwidth:check:level1:level1} ensures the fulfillment of inequalities~\eqref{eqn:bandwidth:check:level1:level2}, according to relations~\eqref{eqn:relation:bandwidth:level-1}.

The partial check approach is more efficient than the direct check approach. When assigning the working route $r_{k}^{0 0}$ to demand $k$, the number of inequalities to check reduces from $|\set{L}(r_{k}^{0 0})| (1 + |\set{E}|^{2})$ to $|\set{L}(r_{k}^{0 0})|$, and when assigning the level-1 route $r_{k}^{e_1 0}$ to demand $k$, the number of inequalities to check decreases from $|\set{L}(r_{k}^{e_1 0})| (1 + |\set{E} \setminus \set{E}(r_{k}^{e_1 0})|)$ to $|\set{L}(r_{k}^{e_1 0})|$.

In the partial check approach, only bandwidth capacity in the scenario being currently planned needs verification. Therefore, when finding a route for demand $k$ in scenario $(e_1, e_2) \in \Omega$ with the labeling algorithm, to determine whether the extension of a label $L_1$ along a logical arc $f_2 \in \set{A}_l$ satisfies bandwidth capacity constraints, we only need to check the inequality $C_{l(f_2)}^{e_1 e_2} + b_k \le B$ for $(e_1, e_2) = (0, 0)$ and the inequality $C_{l(f_2)}^{e_1 e_2} + b_k (1 - y_{k l(f_2)}^{0 0}) \le B$ for $(e_1, e_2) \in \Omega^{1} \cup \Omega^{2}$.

\subsection{Efficient Update of Bandwidth Consumption}
\label{sect:bandwidth:update}

After assigning a route $r_{k}^{e_1 e_2}$ to a demand $k \in \set{K}'$ in scenario $(e_1, e_2) \in \Omega$, the bandwidth consumption on associated lightpaths $\vect{C}_{l}, l \in \set{L}(r_{k}^{e_1 e_2})$ need to be updated for their usage in the check of bandwidth capacity constraints. For this, we first introduce a direct update approach, which is time-consuming, and then present an incremental update approach, which requires significantly fewer update operations and less memory.

The direct update approach stores the bandwidth consumption $C_{l}^{e_1 e_2}$ for each lightpath $l \in \set{L}'$ and each scenario $(e_1, e_2) \in \Omega$. When assigning a route to demand $k \in \set{K}'$ in a scenario, the bandwidth consumption on lightpaths is updated according to Table~\ref{tab:bandwidth:update}, where operations~\eqref{eqn:bandwidth:update:working:failure} are performed for bandwidth reservation constraints, and operations~\eqref{eqn:bandwidth:update:level1:level2} are performed due to consistent routing constraints.
\begin{subequations}
	\begin{eqnarray}
		&& C_{l}^{0 0} \gets C_{l}^{0 0} + b_k, \qquad \forall l \in \set{L}(r_{k}^{0 0}). \label{eqn:bandwidth:update:working:working} \\
		&& C_{l}^{e_1 e_2} \gets C_{l}^{e_1 e_2} + b_k, \qquad \forall l \in \set{L}(r_{k}^{0 0}), (e_1, e_2) \in \Omega^{1} \cup \Omega^{2}. \label{eqn:bandwidth:update:working:failure} \\
		&& C_{l}^{e_1 0} \gets C_{l}^{e_1 0} + b_k (1 - y_{k l}^{0 0}), \qquad \forall l \in \set{L}(r_{k}^{e_1 0}). \label{eqn:bandwidth:update:level1:level1} \\
		&& C_{l}^{e_1 e_2} \gets C_{l}^{e_1 e_2} + b_k (1 - y_{k l}^{0 0}), \qquad \forall l \in \set{L}(r_{k}^{e_1 0}), (e_1, e_2) \in \bar{\Theta}_{k}(e_1). \label{eqn:bandwidth:update:level1:level2} \\
		&& C_{l}^{e_1 e_2} \gets C_{l}^{e_1 e_2} + b_k (1 - y_{k l}^{0 0}), \qquad \forall l \in \set{L}(r_{k}^{e_1 e_2}). \label{eqn:bandwidth:update:level2}
	\end{eqnarray}
\end{subequations}

\begin{table}
\centering
\caption{Operations for updating bandwidth consumption \label{tab:bandwidth:update}}
	{\begin{tabular}{lll}
			\hline
			& Direct Update & Incremental Update \\
			\hline
			Assign $r_{k}^{0 0}$ in scenario $(0, 0)$ & \eqref{eqn:bandwidth:update:working:working} and \eqref{eqn:bandwidth:update:working:failure} & \eqref{eqn:bandwidth:update:working:working} \\
			Assign $r_{k}^{e_1 0}$ in scenario $(e_1, 0) \in \Phi_k$ & \eqref{eqn:bandwidth:update:level1:level1} and \eqref{eqn:bandwidth:update:level1:level2} & \eqref{eqn:bandwidth:update:level1:level1:increment} and \eqref{eqn:bandwidth:update:level1:level2:over-reserve} \\
			Assign $r_{k}^{e_1 e_2}$ in scenario $(e_1, e_2) \in \Psi_k \cup \Theta_k$ & \eqref{eqn:bandwidth:update:level2} & \eqref{eqn:bandwidth:update:level2} \\
			\hline
	\end{tabular}}
\end{table}

The large cardinalities of $\set{E}$ and $\set{E} \setminus \set{E}(r_{k}^{e_1 0})$ in large-scale STG2 instances can lead to a significant number of operations~\eqref{eqn:bandwidth:update:working:failure} and \eqref{eqn:bandwidth:update:level1:level2} ($|\set{L}(r_{k}^{0 0})| |\set{E}|^{2}$ and $|\set{L}(r_{k}^{e_1 0})| |\set{E} \setminus \set{E}(r_{k}^{e_1 0})|$, respectively). However, in Algorithm~\ref{alg:construction:outline}, bandwidth consumption can be updated with the incremental update approach specified in Table~\ref{tab:bandwidth:update}, reducing both memory usage and the number of update operations. Here, $\Delta_{l}^{e_1 e_2} = C_{l}^{e_1 e_2} - C_{l}^{0 0}, l \in \set{L}', (e_1, e_2) \in \Omega^{1} \cup \Omega^{2}$ represents the incremental bandwidth consumption in scenario $(e_1, e_2)$ compared to the working scenario, and $\delta_{l}^{e_1 e_2} = \Delta_{l}^{e_1 0} - \Delta_{l}^{e_1 e_2}, l \in \set{L}', (e_1, e_2) \in \Omega^{2}$ indicates the over-reserved bandwidth in scenario $(e_1, e_2)$ compared to the level-1 scenario $(e_1, 0)$.

Our incremental update approach works as follows. According to Proposition~\ref{prop:relation:bandwidth}, when assigning route $r_{k}^{0 0}$ to demand $k \in \set{K}'$ in the working scenario (lines~\ref{alg:construction:outline:working:begin}--\ref{alg:construction:outline:working:end} of Algorithm~\ref{alg:construction:outline}), we perform operations~\eqref{eqn:bandwidth:update:working:working} and keep $\Delta$ as a zero vector, which is equivalent to performing operations~\eqref{eqn:bandwidth:update:working:working} and \eqref{eqn:bandwidth:update:working:failure}.

By the definition of $\Delta$, we know that operations~\eqref{eqn:bandwidth:update:level1:level1} and \eqref{eqn:bandwidth:update:level1:level2} are equivalent to \eqref{eqn:bandwidth:update:level1:level1:increment} and \eqref{eqn:bandwidth:update:level1:level2:increment}.
\begin{subequations}
	\begin{eqnarray}
		&& \Delta_{l}^{e_1 0} \gets \Delta_{l}^{e_1 0} + b_k (1 - y_{k l}^{0 0}), \qquad \forall l \in \set{L}(r_{k}^{e_1 0}), \label{eqn:bandwidth:update:level1:level1:increment} \\
		&& \Delta_{l}^{e_1 e_2} \gets \Delta_{l}^{e_1 e_2} + b_k (1 - y_{k l}^{0 0}), \qquad \forall l \in \set{L}(r_{k}^{e_1 0}), (e_1, e_2) \in \bar{\Theta}_{k}(e_1). \label{eqn:bandwidth:update:level1:level2:increment}
	\end{eqnarray}
\end{subequations}
Additionally, by the definition of $\delta$, performing operations~\eqref{eqn:bandwidth:update:level1:level1:increment} and \eqref{eqn:bandwidth:update:level1:level2:increment} is equivalent to performing \eqref{eqn:bandwidth:update:level1:level1:increment} and \eqref{eqn:bandwidth:update:level1:level2:over-reserve}.
\begin{equation}
	\delta_{l}^{e_1 e_2} \gets \delta_{l}^{e_1 e_2} + b_k (1 - y_{k l}^{0 0}), \qquad \forall l \in \set{L}(r_{k}^{e_1 0}), (e_1, e_2) \in \Theta_{k}(e_1). \label{eqn:bandwidth:update:level1:level2:over-reserve}
\end{equation}

The number of operations~\eqref{eqn:bandwidth:update:level1:level2:over-reserve}, i.e., $|\set{L}(r_{k}^{e_1 0})| |\Theta_{k}(e_1)| = |\set{L}(r_{k}^{e_1 0})| |\set{E}(r_{k}^{e_1 0})|$, can be significantly smaller than that of operations~\eqref{eqn:bandwidth:update:level1:level2:increment} (i.e., $|\set{L}(r_{k}^{e_1 0})| |\set{E} \setminus \set{E}(r_{k}^{e_1 0})|$) in large-scale STG2 instances, since we often have $|\set{E}| \gg |\set{E}(r_{k}^{e_1 0})|$ in these instances. Therefore, when assigning route $r_{k}^{e_1 0}$ to demand $k \in \set{K}'$ in a level-1 scenario $(e_1, 0) \in \Phi_k$ (lines~\ref{alg:construction:outline:level-1:begin}--\ref{alg:construction:outline:level-1:end} of Algorithm~\ref{alg:construction:outline}), the query of bandwidth consumption $C_{l}^{e_1 0}$ is given by $C_{l}^{0 0} + \Delta_{l}^{e_1 0}$ for $l \in \set{L}'$, while bandwidth consumption update is performed by operations~\eqref{eqn:bandwidth:update:level1:level1:increment} and \eqref{eqn:bandwidth:update:level1:level2:over-reserve}.

At the beginning of planning a level-2 scenario $(e_1, e_2) \in \Omega^{2}$, we initialize bandwidth consumption $C_{l}^{e_1 e_2}$ for $l \in \set{L}'$ with operation~\eqref{eqn:bandwidth:level-2:init}. When assigning route $r_{k}^{e_1 e_2}$ to demand $k \in \set{K}'$ in scenario $(e_1, e_2) \in \Psi_k \cup \Theta_k$, the bandwidth consumption is updated according to operations~\eqref{eqn:bandwidth:update:level2}. In Algorithm~\ref{alg:construction:outline}, level-2 scenarios are planned sequentially (lines~\ref{alg:construction:outline:level-2:begin}--\ref{alg:construction:outline:level-2:end}), and planning a level-2 scenario imposes no constraint on another level-2 scenario. Therefore, at the end of planning scenario $(e_1, e_2) \in \Omega^{2}$, we release the memory for storing bandwidth consumption $C_{l}^{e_1 e_2}, \forall l \in \set{L}'$.
\begin{equation}
	\label{eqn:bandwidth:level-2:init}
	C_{l}^{e_1 e_2} \gets C_{l}^{0 0} + \Delta_{l}^{e_1 0} - \delta_{l}^{e_1 e_2}.
\end{equation}

Compared with the direct update approach, the incremental update approach requires fewer operations. Specifically, when assigning the working route $r_{k}^{0 0}$ to demand $k$, the number of update operations decreases from $|\set{L}(r_{k}^{0 0})| (1 + |\set{E}|^{2})$ to $|\set{L}(r_{k}^{0 0})|$, and when assigning the level-1 route $r_{k}^{e_1 0}$ to demand $k$, the number of update operations decreases from $|\set{L}(r_{k}^{e_1 0})| (1 + |\set{E} \setminus \set{E}(r_{k}^{e_1 0})|)$ to $|\set{L}(r_{k}^{e_1 0})| (1 + |\set{E}(r_{k}^{e_1 0})|)$, where $|\set{E} \setminus \set{E}(r_{k}^{e_1 0})|$ is often much larger than $|\set{E}(r_{k}^{e_1 0})|$ in large-scale STG2 instances.

Furthermore, for each lightpath $l \in \set{L}'$, the incremental update approach only maintains $C_{l}^{0 0}$, $\Delta_{l}^{e_1 0}, \forall (e_1, 0) \in \bigcup_{k \in \set{K}'} \Phi_{k}$, and $\delta_{l}^{e_1 e_2}, \forall (e_1, e_2) \in \bigcup_{k \in \set{K}'} \Theta_{k}$, along with an additional scalar $C_{l}^{e_1 e_2}$ when $(e_1, e_2)$ is the level-2 scenario being currently planned. These data are often much sparser than that maintained in the direct update approach (i.e., $C_{l}^{e_1 e_2}, \forall (e_1, e_2) \in \Omega$), indicating that the incremental update approach also requires less memory.

\section{Details about Parallel Enhancement}

\subsection{Planning with Slave Threads}
\label{sect:parallel:slave}

Multiple slave threads are used to plan level-2 scenarios in parallel without generating new lightpaths, while within a single thread scenarios are planned sequentially. The $t$-th thread executes the procedure shown in Algorithm~\ref{alg:parallel:slave} to plan level-2 scenarios $\Gamma_t$. The successfully planned scenarios are denoted as $\bar{\Gamma}_t$, while $\tilde{\Gamma}_t$ represent the set of blocked scenarios (failed to be planned due to insufficient lightpaths). Blocked scenarios $\tilde{\Gamma}_t$ are then planned in the master thread, where new lightpaths can be established.

Algorithm~\ref{alg:parallel:slave} proceeds as follows. The values to be returned are initialized in line~\ref{alg:parallel:slave:init:output}. In line~\ref{alg:parallel:slave:init:stop}, the indicator to stop the thread $flag_{stop}$ and the number of routes generated in this thread $n_{R}$ are initialized. At the beginning of planning a level-2 scenario $(e_1, e_2)$ in line~\ref{alg:parallel:slave:init:bandwidth}, bandwidth consumption on lightpaths $C_{l}^{e_1 e_2}, l \in \set{L}'$ is initialized with \eqref{eqn:bandwidth:level-2:init}. For each unprotected demand in scenario $(e_1, e_2)$, we attempt to find the level-2 route by evaluating generated routes (line~\ref{alg:parallel:slave:reuse}) or using the labeling algorithm without establishing new lightpaths (line~\ref{alg:parallel:slave:labeling}). If both are failed, the scenario should be planned in the master thread, and we record the bandwidth consumption on lightpaths in this scenario (line~\ref{alg:parallel:slave:record-failed-scenario}).

Stopping criteria are evaluated in lines~\ref{alg:parallel:slave:check-stop:fail}, \ref{alg:parallel:slave:check-stop:route}, and \ref{alg:parallel:slave:check-stop:scenario}. The indicator of whether there is an idle slave thread $flag_{idle}$, the total number of scenarios that failed to be planned in all slave threads $n_{fail}$, the total number of routes newly generated in all slave threads $n_{route}$, are modifiable data shared by all slave threads. Hyperparameters $NF$ and $NR$ represent thresholds on $n_{fail}$ and $n_{route}$, respectively. The result of evaluating stopping criteria is $true$ if and only if at least one of the following conditions is satisfied: (i) $flag_{idle} = true$; (ii) $n_{fail} \ge NF$; (iii) $n_{route} \ge NR$.

The motivation for stopping criterion (i) is to reduce thread idleness: When a thread has already stopped, stopping all slave threads earlier allows an earlier start of a new ``master-slave'' iteration of Algorithm~\ref{alg:parallel:outline} and an earlier reactivation of stopped threads. Stopping criteria (ii) and (iii) aim to control the number of ``master-slave'' iterations. Smaller $NF$ and $NR$ lead to more iterations, increasing the overhead of starting and stopping threads. Conversely, a larger $NF$ results in fewer iterations, but tends to increase the number of scenarios being planned by the master thread, reducing parallelization. A larger $NR$ results in less frequent synchronization of the route pool $\set{R}'$ and more frequent calls to the labeling algorithm.

Lines~\ref{alg:parallel:slave:check-check:fail}, \ref{alg:parallel:slave:check-check:route}, and \ref{alg:parallel:slave:check-check:scenario} determine when to evaluate stopping criteria. The hyperparameters $CS$ and $CR$ control this frequency. A higher frequency can lead to race conditions, i.e., multiple threads access and modify shared data ($flag_{idle}, n_{fail}, n_{route}$) concurrently. Whereas a lower frequency can cause idleness, as already stopped slave threads remain idle until all slave threads stop.
{
\begin{algorithm}
\caption{Slave Plan}
\label{alg:parallel:slave}
\vskip6pt
{\small
\begin{algorithmic}[1]

	\Statex \algorithmicrequire \text{ } $\Gamma_t$, hyperparameters $NF, NR, CS, CR$, and modifiable data shared by all slave threads $flag_{idle}, n_{fail}, n_{route}$
	\Statex \algorithmicensure \text{ } $\bar{\Gamma}_t, \tilde{\Gamma}_t, \tilde{\vect{C}}(t)$
	\State $\bar{\Gamma}_t \gets \emptyset, \tilde{\Gamma}_t \gets \emptyset, \tilde{\vect{C}}(t) \gets \emptyset, \set{R}^{\prime}(t) \gets \emptyset$ \label{alg:parallel:slave:init:output}
	\State $flag_{stop} \gets false, n_{R} \gets 0$ \label{alg:parallel:slave:init:stop}
	\ForAll{$(e_1, e_2) \in \Gamma_t$ \algand $flag_{stop} = false$}
		\State $flag_{fail} \gets false$ $\qquad \qquad \quad$ (Whether fails to plan some demand in the scenario)
		\State Initialize bandwidth consumption on lightpaths $C_{l}^{e_1 e_2}, l \in \set{L}'$ with \eqref{eqn:bandwidth:level-2:init} \label{alg:parallel:slave:init:bandwidth}
		\ForAll{$k \in \set{K}'$ \algand demand $k$ is not protected in scenario $(e_1, e_2)$}
			\If{found a generated route $r \in \set{R}'$ for $k$ in $(e_1, e_2)$ with Algorithm~\ref{alg:reuse-route}} \label{alg:parallel:slave:reuse}
			\State $r_{k}^{e_1 e_2} \gets r$; Update bandwidth consumption with \eqref{eqn:bandwidth:update:level2}
			\ElsIf{found a new route $r$ for $k$ in $(e_1, e_2)$ with Algorithm~\ref{alg:labeling}} \label{alg:parallel:slave:labeling}
				\State $r_{k}^{e_1 e_2} \gets r$; Update bandwidth consumption with \eqref{eqn:bandwidth:update:level2}
				\State $\set{R}^{\prime}(t) \gets \set{R}^{\prime}(t) \cup \{r\}$
			\Else
				\State $flag_{fail} \gets true$
			\EndIf
		\EndFor
		\If{$flag_{fail} = false$}
			\State $\bar{\Gamma}_t \gets \bar{\Gamma}_t \cup \{(e_1, e_2)\}$ $\qquad \qquad \qquad \quad$ (All demands are planned in the scenario)
		\Else
			\State $\tilde{\Gamma}_t \gets \tilde{\Gamma}_t \cup \{(e_1, e_2)\}$ $\qquad \qquad \qquad \qquad \qquad \qquad \qquad$ (The scenario is blocked)
			\State $\tilde{\vect{C}}(t) \gets \tilde{\vect{C}}(t) \cup \{C_{l}^{e_1 e_2}: l \in \set{L}'\}$ \label{alg:parallel:slave:record-failed-scenario} $\qquad \qquad \qquad$ (Store bandwidth consumption)
		\EndIf
		\If{$flag_{fail} = true$} \label{alg:parallel:slave:check-check:fail}
			\State $n_{fail} \gets n_{fail} + 1$
			\State $flag_{stop} \gets$ \Call{StopCriteria}{$NF, NR, flag_{idle}, n_{fail}, n_{route}$} \label{alg:parallel:slave:check-stop:fail}
		\ElsIf{$|\set{R}^{\prime}(t)| - n_{R} \ge CR$} \label{alg:parallel:slave:check-check:route}
			\State $n_{route} \gets n_{route} + |\set{R}^{\prime}(t)| - n_{R}$
			\State $n_{R} \gets |\set{R}^{\prime}(t)|$
			\State $flag_{stop} \gets$ \Call{StopCriteria}{$NF, NR, flag_{idle}, n_{fail}, n_{route}$} \label{alg:parallel:slave:check-stop:route}
		\ElsIf{$(|\bar{\Gamma}_t| + |\tilde{\Gamma}_t|) \% CS = 0$} \label{alg:parallel:slave:check-check:scenario}
			\State $flag_{stop} \gets$ \Call{StopCriteria}{$NF, NR, flag_{idle}, n_{fail}, n_{route}$} \label{alg:parallel:slave:check-stop:scenario}
		\EndIf
	\EndFor
	\State $flag_{idle} \gets true$
	\State \Return $\bar{\Gamma}_t, \tilde{\Gamma}_t, \tilde{\vect{C}}(t)$

\end{algorithmic}
}
\end{algorithm}
}
\subsection{Planning with Master Thread}
\label{sect:parallel:master}

By applying the procedure shown in Algorithm~\ref{alg:parallel:master}, the master thread plans level-2 scenarios that slave threads cannot complete. The main differences between planning with the master thread and slave threads are: (i) The bandwidth consumption $C_{l}^{e_1 e_2}$ on lightpath $l \in \set{L}'$ in scenario $(e_1, e_2)$ is initialized with $\tilde{C}_{l}^{e_1 e_2}$ (line~\ref{alg:parallel:master:init:bandwidth}) instead of using \eqref{eqn:bandwidth:level-2:init}; (ii) New lightpaths can be established when finding routes with the labeling algorithm (line~\ref{alg:parallel:master:labeling}); (iii) If the master thread cannot find a level-2 route to protect a demand in a scenario, the hierarchical constructive heuristic fails to find a feasible solution to the STG2 (line~\ref{alg:parallel:master:fail}).

\begin{algorithm}
\caption{Master Plan}
\label{alg:parallel:master}
\vskip6pt
{\small
\begin{algorithmic}[1]

	\Statex \algorithmicrequire \text{ } blocked scenarios $\tilde{\Gamma}$, bandwidth consumption in blocked scenarios $\tilde{\vect{C}}$
	\Statex \algorithmicensure \text{ } Routes $r_{k}^{e_1 e_2}, k \in \set{K}, (e_1, e_2) \in \tilde{\Gamma}$
	\ForAll{$(e_1, e_2) \in \tilde{\Gamma}$}
		\State Initialize bandwidth consumption on lightpaths $C_{l}^{e_1 e_2} = \tilde{C}_{l}^{e_1 e_2}, l \in \set{L}'$ \label{alg:parallel:master:init:bandwidth}
		\ForAll{$k \in \set{K}'$ \algand demand $k$ is not protected in scenario $(e_1, e_2)$}
			\If{found a generated route $r \in \set{R}'$ for $k$ in $(e_1, e_2)$ with Algorithm~\ref{alg:reuse-route}}
				\State $r_{k}^{e_1 e_2} \gets r$; Update bandwidth consumption with \eqref{eqn:bandwidth:update:level2}
			\ElsIf{found a new route $r$ for $k$ in $(e_1, e_2)$ with Algorithm~\ref{alg:labeling}} \label{alg:parallel:master:labeling}
				\State $r_{k}^{e_1 e_2} \gets r$; Update bandwidth consumption with \eqref{eqn:bandwidth:update:level2}
				\State Insert new lightpaths into $\set{L}'$, if any.
			\Else
				\State Fails to find a feasible solution; Terminate. \label{alg:parallel:master:fail}
			\EndIf
		\EndFor
	\EndFor
\end{algorithmic}
}
\end{algorithm}

\section{Computing Time for Planning Level-2 Scenarios in Parallel}
\label{sect:experiments:parallel}

Table~\ref{tab:experiments:effectiveness:parallel} displays the computing time for planning level-2 scenarios with different numbers of threads, as shown in column ``Time2''.

\begin{table}[H]
\centering
\caption{Computational results of parallel algorithms with different number of threads \label{tab:experiments:effectiveness:parallel}}
{\begin{tabular}{rrrr|rr|rr|rr}
\hline
\multicolumn{4}{c|}{Instance} & \multicolumn{2}{c|}{CH-1} & \multicolumn{2}{c|}{CH-4} & \multicolumn{2}{c}{CH-8} \\
\hline
$|\set{N}|$ & $|\set{E}|$ & $|\set{K}|$ & $\bar{d}$ (km) & Obj & Time2 (s) & Obj & Time2 (s) & Obj & Time2 (s) \\
\hline
1000 & 1928 & 10020 & 360 & 5959 & 1674.3 & 5862 & 641.5 & 5882 & 344.6 \\
1000 & 1928 & 10020 & 600 & 5883 & 1671.6 & 5863 & 646.7 & 5870 & 345.7 \\
1000 & 1928 & 10020 & 900 & 5179 & 1252.8 & 5176 & 515.8 & 5131 & 261.4 \\
1000 & 1928 & 10020 & 1500 & 4763 & 1056.1 & 4805 & 415.5 & 4792 & 217.6 \\
1200 & 2365 & 10390 & 600 & 6675 & 2637.1 & 6589 & 1005.2 & 6605 & 498.1 \\
1200 & 2365 & 10390 & 900 & 5903 & 1972.6 & 5856 & 803.3 & 5874 & 398.8 \\
1200 & 2365 & 10390 & 1500 & 5513 & 1602 & 5486 & 659 & 5514 & 334 \\
1600 & 3091 & 10070 & 600 & 6700 & 4679 & 6673 & 1624.8 & 6658 & 748.6 \\
1600 & 3091 & 10070 & 900 & 5804 & 5204.2 & 5750 & 996.4 & 5779 & 584.4 \\
2000 & 3083 & 12980 & 600 & 8543 & 14262.6 & 8496 & 2626.3 & 8540 & 1552.7 \\
2200 & 3052 & 11054 & 600 & 8925 & 17432.4 & 8839 & 3956.4 & 8882 & 1880.8 \\
2200 & 3052 & 11054 & 900 & 6858 & 10437 & 6850 & 1925.6 & 6856 & 1153.9 \\
2200 & 3052 & 11054 & 1500 & 5591 & 4754.8 & 5545 & 1418.4 & 5541 & 831.9 \\
2600 & 3027 & 10930 & 600 & 8924 & 9547.2 & 8822 & 1908.2 & 8908 & 1169.8 \\
2600 & 3027 & 10930 & 900 & 6588 & 5258.2 & 6559 & 1228.9 & 6573 & 742.4 \\
2600 & 3027 & 10930 & 1500 & 4990 & 3571.7 & 4953 & 911.4 & 4942 & 539.8 \\
\hline
\multicolumn{4}{c|}{Average $TR_1$} & & 1.0 & & 3.6 & & 6.5 \\
\hline
\end{tabular}}
\end{table}

\section{Computational Results for Evaluating Algorithm Enhancements}
\label{app:sec:enhance}
In Table \ref{tab:experiments:effectiveness:aggregation} and Table~\ref{tab:experiments:effectiveness:routing}, compared to CH, $Imp_C$ represents the improvement in objective value, and $TR_C$ is the ratio of computing times. Table \ref{tab:experiments:effectiveness:aggregation} evaluates the effectiveness of demand aggregation. The number of aggregated demands is shown in column $|\set{K}'|$, which is approximately one-tenth of the original number. Table~\ref{tab:experiments:effectiveness:routing} evaluate the performance of different routing techniques. For CH-R-DP and CH-R, due to prohibitive computing times in the first seven instances, we do not test larger-scale instances. 

\begin{table}[H]
\centering
\caption{Computational results of disabling/allowing demand aggregation \label{tab:experiments:effectiveness:aggregation}}
{\begin{tabular}{rrrr|rrrr|rrrr}
\hline
\multicolumn{4}{c|}{Instance} & \multicolumn{4}{c|}{CH-D} & \multicolumn{3}{c}{CH} \\
\hline
$|\set{N}|$ & $|\set{E}|$ & $|\set{K}|$ & $\bar{d}$ (km) & Obj & $Imp_C$ & Time (s) & $TR_{C}$ & Obj  & Time (s) & $|\set{K}'|$ \\
\hline
1000 & 1928 & 10020 & 360 & 5912 & -0.5\% & 1526.8 & 3.9 & 5882 & 393.6 & 1002 \\
1000 & 1928 & 10020 & 600 & 5914 & -0.7\% & 1490 & 3.8 & 5870 & 394.1 & 1002 \\
1000 & 1928 & 10020 & 900 & 5152 & -0.4\% & 1249.9 & 4.1 & 5131 & 305.2 & 1002 \\
1000 & 1928 & 10020 & 1500 & 4749 & 0.9\% & 1078.6 & 4.1 & 4792 & 262.4 & 1002 \\
1200 & 2365 & 10390 & 600 & 6705 & -1.5\% & 2024.3 & 3.6 & 6605 & 564.7 & 1039 \\
1200 & 2365 & 10390 & 900 & 5877 & -0.1\% & 1692.4 & 3.7 & 5874 & 458 & 1039 \\
1200 & 2365 & 10390 & 1500 & 5418 & 1.7\% & 1450.6 & 3.7 & 5514 & 394.8 & 1039 \\
1600 & 3091 & 10070 & 600 & 6658 & 0.0\% & 2911.2 & 3.5 & 6658 & 837.4 & 1007 \\
1600 & 3091 & 10070 & 900 & 5630 & 2.6\% & 2426.3 & 3.7 & 5779 & 663.8 & 1007 \\
2000 & 3083 & 12980 & 600 & 9401 & -10.1\% & 6078.9 & 3.6 & 8540 & 1695 & 1298 \\
2200 & 3052 & 11054 & 600 & 10638 & -19.8\% & 9265 & 4.6 & 8882 & 2031.5 & 1106 \\
2200 & 3052 & 11054 & 900 & 7734 & -12.8\% & 5106.7 & 4 & 6856 & 1284.1 & 1106 \\
2200 & 3052 & 11054 & 1500 & 5839 & -5.4\% & 3686.1 & 3.9 & 5541 & 949.1 & 1106 \\
2600 & 3027 & 10930 & 600 & 11687 & -31.2\% & 9670.7 & 7.7 & 8908 & 1254.4 & 1028 \\
2600 & 3027 & 10930 & 900 & 8240 & -25.4\% & 4848.2 & 5.9 & 6573 & 817.6 & 1028 \\
2600 & 3027 & 10930 & 1500 & 5855 & -18.5\% & 3515.6 & 5.8 & 4942 & 607.2 & 1028 \\
\hline
\multicolumn{4}{c|}{Average $Imp_C$ and $TR_C$} & & -7.6\% & & 4.3 & & & \\
\hline
\end{tabular}}
\end{table}


\begin{table}[H]
\begin{small}
\centering
\caption{Computational results of adopting different routing techniques \label{tab:experiments:effectiveness:routing}}
{\begin{tabular}{rrrr|rr|rr|rr|rr}
\hline
\multicolumn{4}{c|}{Instance} & \multicolumn{2}{c|}{CH-R-DP} & \multicolumn{2}{c|}{CH-R} & \multicolumn{2}{c|}{CH-B} & \multicolumn{2}{c}{CH} \\
\hline
$|\set{N}|$ & $|\set{E}|$ & $|\set{K}|$ & $\bar{d}$ (km) & Obj & Time (s) & Obj & Time (s) & Obj & Time (s) & Obj & Time (s) \\
\hline
1000 & 1928 & 10020 & 360  & 6038 & 18089.5 & 6108 & 14603.6 & 5865 & 601.1  & 5882 & 393.6  \\
1000 & 1928 & 10020 & 600  & 6082 & 17702.2 & 6039 & 16175.2 & 5918 & 606.9  & 5870 & 394.1  \\
1000 & 1928 & 10020 & 900  & -    & -       & 5339 & 12102.8 & 5199 & 492.0  & 5131 & 305.2  \\
1000 & 1928 & 10020 & 1500 & -    & -       & 4903 & 9062.8  & 4757 & 462.3  & 4792 & 262.4  \\
1200 & 2365 & 10390 & 600  & 6762 & 32498.2 & 6752 & 25055.3 & 6640 & 934.1  & 6605 & 564.7  \\
1200 & 2365 & 10390 & 900  & -    & -       & 5995 & 20599.8 & 5855 & 743.8  & 5874 & 458.0  \\
1200 & 2365 & 10390 & 1500 & -    & -       & 5635 & 13360.9 & 5556 & 641.1  & 5514 & 394.8  \\
1600 & 3091 & 10070 & 600  & *    & *       & *    & *       & 6700 & 1223.7 & 6658 & 837.4  \\
1600 & 3091 & 10070 & 900  & *    & *       & *    & *       & 5709 & 974.1  & 5779 & 663.8  \\
2000 & 3083 & 12980 & 600  & *    & *       & *    & *       & 8541 & 3149.8 & 8540 & 1695.0 \\
2200 & 3052 & 11054 & 600  & *    & *       & *    & *       & 8961 & 5693.6 & 8882 & 2031.5 \\
2200 & 3052 & 11054 & 900  & *    & *       & *    & *       & 6806 & 3684.7 & 6856 & 1284.1 \\
2200 & 3052 & 11054 & 1500 & *    & *       & *    & *       & 5557 & 2109.9 & 5541 & 949.1  \\
2600 & 3027 & 10930 & 600  & *    & *       & *    & *       & 8829 & 3467.2 & 8908 & 1254.4 \\
2600 & 3027 & 10930 & 900  & *    & *       & *    & *       & 6574 & 2150.5 & 6573 & 817.6  \\
2600 & 3027 & 10930 & 1500 & *    & *       & *    & *       & 4947 & 1503.0 & 4942 & 607.2  \\
\hline
\multicolumn{4}{c|}{Average $Imp_C$ and $TR_C$} & -2.8\% & 49.5 & -2.7\% & 39.4 & -0.1\% & 2.0 & & \\
\hline
\end{tabular}}
{Notes. -: Out of memory of 16 GB; *: Not tested, as the algorithm variant cannot solve smaller-scale instances (with 1,000 nodes) within the memory and time limits.}
\end{small}
\end{table}

\end{appendices}

\end{document}